\documentclass [10pt,reqno]{amsart}
\usepackage {color,graphicx}
\usepackage{texdraw}
\usepackage{amsmath}
\usepackage{amssymb} 
\usepackage [latin1]{inputenc}
\usepackage[all]{xy}
\usepackage{epic}
\usepackage{amsthm}
\newtheorem{theorem}{Theorem}
\newtheorem{lemma}[theorem]{Lemma}
\newtheorem{proposition}[theorem]{Proposition}
\newtheorem{remark}[theorem]{Remark}
\newtheorem{example}[theorem]{Example}
\newtheorem{definition}[theorem]{Definition}
\newtheorem{notation}[theorem]{Notation}
\newcommand{\C}{\mathbb C}
\newcommand{\R}{\mathbb R}
\newcommand{\Q}{\mathbb Q}
\newcommand{\Z}{\mathbb Z}
\newcommand{\K}{\mathbb K}
\newcommand{\TP}{\mathbb P}
\newcommand{\bx}{{\boldsymbol x}}
\newcommand{\bp}{{\boldsymbol p}}

\DeclareMathOperator{\Gl}{Gl}
\DeclareMathOperator{\Ker}{Ker}
\DeclareMathOperator{\Sing}{Sing}
\DeclareMathOperator{\Tor}{Tor}
\DeclareMathOperator {\rowspace}{rowspace}
\DeclareMathOperator {\Span}{span}
\DeclareMathOperator {\val}{val}

\DeclareMathOperator {\trop}{trop}
\DeclareMathOperator {\Trop}{Trop}
\DeclareMathOperator {\conv}{conv}

\newcommand{\longer}[1]{}

\title [Tropical surface singularities]{Tropical surface singularities}

\author {Hannah Markwig}
\address {Hannah Markwig, Universit\"at des Saarlandes, Fachrichtung Mathematik, Postfach 151150, 66041 Saarbr\"ucken, Germany }
\email {hannah@math.uni-sb.de}

\author{Thomas Markwig}
\address{Thomas Markwig, TU Kaiserslautern, Erwin-Schr\"odinger-Stra\ss{}e, 67653 Kaisers-lautern, Germany }
\email {keilen@mathematik.uni-kl.de}

\author{Eugenii Shustin}
\address{Eugenii Shustin, School of Mathematical Sciences, Tel Aviv University, Ramat Aviv, Tel Aviv 69978, Israel}
\email {shustin@post.tau.ac.il}


\begin {document}

   \maketitle

   \begin{abstract}
In this paper, we study tropicalisations of singular surfaces in
toric threefolds. We completely classify singular tropical surfaces
of maximal-dimensional geometric type, show that they can generically have
only finitely many singular points, and describe all possible
locations of singular points. More precisely, we show that singular
points must be either vertices, or generalized midpoints and
barycenters of certain faces of singular tropical surfaces, and, in
some case, there may be additional metric restrictions to faces of
singular tropical surfaces.
   \end{abstract}

\section{Introduction}

This paper studies singularities of tropical surfaces in $\R^3$. The
question what the analogue of a singularity in the tropical world
should be is quite natural to ask and has consequently interested
several authors recently (\cite{DFS05}, \cite{DT10}, \cite{MMS09}).
The fact that this question is hard to answer in general makes it
even more intriguing. We define a point $p$ in a tropical surface
$S$ to be singular if there is an algebraic surface $\tilde{S}$,
defined over the Puiseux-series with coefficients in $\C$, whose tropicalisation is $S$ and
which is singular at a point $\tilde{p}\in \tilde{S}$ that
tropicalises to $p$. Given a non-degenerate lattice polytope
$\Delta\in \R^3$, consider the family $\Sing(\Delta)\subset
\TP^{\#(\Delta\cap \Z^3)-1}$ of singular hypersurfaces in the toric
threefold defined by $\Delta$ whose defining equations have Newton
polytope $\Delta$. We assume that $\Delta$ is non-defective, i.e.\
that $\Sing(\Delta)$ is a hypersurface in $\TP^{\#(\Delta\cap
\Z^3)-1}$, defined by a polynomial which is then called the
discriminant of $\Delta$. The tropicalisation $\Trop(\Sing(\Delta))$
of $\Sing(\Delta)$ has been studied in \cite{DFS05} and is called
the tropical discriminant. While a general member of $\Sing(\Delta)$
has exactly one singular point, namely a node, an analogous
statement is not true in tropical geometry. The reason is that for a
given singular tropical surface, there can be several singular
tropical surfaces tropicalising to it, but such that the respective
singular points tropicalise to different points in the tropical
surface. Consequently, there are also general tropical surfaces with
infinitely many singularities. The subset of singular points of a
tropical surface does not seem to have any nice structure,
in particular it is not a tropical subvariety. Examples 4.5 and 4.3
of \cite{DT10} show tropical curves with infinitely many resp. two
singular points. 
We concentrate on singular tropical surfaces of
maximal-dimensional geometric type (see Definition \ref{def-mdg-type} in Subsection \ref{subsec-hypersurf} for a precise description). These are the
singular tropical surfaces whose parameter space is of the maximal possible dimension equal to $\#(\Delta\cap\Z^3)-2$, which, in particular, equals the dimension of the parameter space of singular algebraic surfaces with Newton polygon $\Delta$. Specifically, such tropical surfaces have only finitely many singular points. We completely classify these singular tropical surfaces and describe possible locations of singular points.

Our study is closely related to \cite{DT10}, which deals with
singular tropical hypersurfaces of any dimension. There, a more
algebraic point of view is taken however: the main result is the
description of tropical singular points in terms of Euler
derivatives, i.e.\ tropical equations are given which a point must
satisfy to be singular. We concentrate more on the geometry of
singular tropical surfaces.

Our paper can be viewed as a sequel to \cite{MMS09}, where we studied
tropical plane curves with a singular point. The main result of
\cite{MMS09} is the classification of singular tropical curves of
maximal-dimensional geometric type. A singular point of a tropical curve of maximal-dimensional geometric
type is either a
``crossing'' of two edges, or a three-valent vertex of multiplicity $3$,
or it is a point on an edge $e$ of weight two which has equal distance
to the two vertices of $e$ (or which satisfies a similar metric
condition respectively).
To derive this result, we used the following methods: we considered the
family of algebraic curves in a toric surface with a singularity in a
fixed point. This family is defined by linear equations, and so its
tropicalisation is a Bergman fan which can be described in terms of
weight classes of flags of flats of the corresponding matroid
(\cite{FS05}, \cite{AK06}). We studied the possible weight classes and
classified the corresponding tropical curves.
Fixing a different point in the torus yields a shift of the Bergman fan
 (see Remark 3.1 and 3.2 of \cite{MMS09}).

Here, we apply the same methods to the family of algebraic surfaces
in a toric threefold with a singularity in a fixed point. While the
basic ideas we use are the same as in \cite{MMS09}, the
classification becomes much more complicated and we have to
establish and use various facts about lattice polytopes. Also, we
concentrate purely on tropical surfaces with only finitely many
singularities (contrary to our classification in the curve case in
\cite{MMS09}). Our main result is the classification in Theorem
\ref{thm2} below. Such a classification is not possible in higher
dimensions (see Remark \ref{rem-discrim}). Theorem \ref{thm1} tells
us for which tropical surfaces there are only finitely many
singularities. For more details and notation, see Section
\ref{sec-class}.

\begin{theorem}\label{thm1}
  Let $\Delta\subset\R^n$ be a non-degenerate convex lattice polytope
  and denote by $\mathcal{A}=\Delta\cap\Z^n$ the lattice points of $\Delta$.
 Let $F_u(\bx)=\max_{m\in\mathcal{A}}\{u_{m}+m\cdot\bx\}$,
  $\bx\in\R^n$,
  define a generic (see Definition \ref{def-generic1}) singular tropical
  hypersurface $S$. Assume the dual marked subdivision corresponds to a
  cone of codimension $c$ in the secondary fan.
  Then the set of singular points in $S$ is a union of finitely many
  polyhedra of dimension $c-1$.
\end{theorem}

In the following classification below, we thus want to restrict to the
case $c=1$ of generic tropical surfaces $S$ whose dual marked subdivision
corresponds to a cone of codimension $1$ in the secondary fan. (Not all cones of codimension $1$ in the secondary fan correspond to singular tropical surfaces. We do not give a complete classification but restrict to cones of maximal-dimensional geometric type.)
It follows that the dual marked subdivision contains a unique circuit
and that every marked polytope in the subdivision which does not
contain the circuit is a simplex (see Remark \ref{rem-circuit}).
We can conclude from Lemma 3.1 of \cite{DT10} that every
singular point of $S$ is contained in the cell of $S$ dual to the
circuit.

In addition, we make the assumption that the
tropical surface is of maximal-dimensional geometric type
(see Definition \ref{def-mdg-type} in Subsection \ref{subsec-hypersurf}).
In this case, the singular tropical surface uniquely defines
a codimension one cone of the secondary fan, and, in the dual marked subdivision, all lattice points of $\Delta$ are marked. Our main result is
a complete classification of such singular tropical surfaces and
of possible locations of their singular points. 

Notice that some codimension 1 cones of the secondary fan do not appear in our
classification: these correspond to singular tropical surfaces which are not of maximal-dimensional geometric type. In this case the cone cannot be uniquely restored out of the tropical surface, and the singular locus has positive dimension.

In what follows we will usually consider polytopes only up to integral
unimodular affine transformations which we refer to as {\em IUA-equivalence}.

\begin{theorem}\label{thm2}
  Let $F_u=\max_{(i,j,k)\in\mathcal{A}}\{u_{(i,j,k)}+ix+jy+kz\}$
  define a singular tropical surface $S$.
  We assume that $S$ is generic (see Definition \ref{def-generic1}) and dual to a
  marked subdivision $T=\{(Q_1,\mathcal{A}_1),\ldots,(Q_k,\mathcal{A}_k)\}$ (see Subsection
  \ref{subsec-hypersurf}) of maximal-dimensional geometric type.
  Assume the dual subdivision corresponds to a cone of codimension $1$ in the secondary fan.
  Then every marked polytope $(Q_i,\mathcal{A}_i)$ in $T$ which does
  not contain the circuit is a simplex, and $S$ contains only finitely many
  singular points.
Their possible locations and dual polytopes, classified up to
IUA-equivalence, are as follows:
  \begin{figure}[h]
    \centering
\begin{picture}(0,0)%
\includegraphics{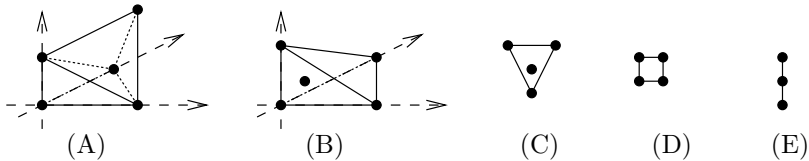}%
\end{picture}%
\setlength{\unitlength}{3947sp}%
\begingroup\makeatletter\ifx\SetFigFont\undefined%
\gdef\SetFigFont#1#2#3#4#5{%
  \reset@font\fontsize{#1}{#2pt}%
  \fontfamily{#3}\fontseries{#4}\fontshape{#5}%
  \selectfont}%
\fi\endgroup%
\begin{picture}(5042,992)(3064,-2966)
\put(7876,-2911){\makebox(0,0)[lb]{\smash{{\SetFigFont{10}{12.0}{\familydefault}{\mddefault}{\updefault}{\color[rgb]{0,0,0}(E)}%
}}}}
\put(3451,-2911){\makebox(0,0)[lb]{\smash{{\SetFigFont{10}{12.0}{\familydefault}{\mddefault}{\updefault}{\color[rgb]{0,0,0}(A)}%
}}}}
\put(4951,-2911){\makebox(0,0)[lb]{\smash{{\SetFigFont{10}{12.0}{\familydefault}{\mddefault}{\updefault}{\color[rgb]{0,0,0}(B)}%
}}}}
\put(6301,-2911){\makebox(0,0)[lb]{\smash{{\SetFigFont{10}{12.0}{\familydefault}{\mddefault}{\updefault}{\color[rgb]{0,0,0}(C)}%
}}}}
\put(7126,-2911){\makebox(0,0)[lb]{\smash{{\SetFigFont{10}{12.0}{\familydefault}{\mddefault}{\updefault}{\color[rgb]{0,0,0}(D)}%
}}}}
\end{picture}%
    \caption{The possible circuits.}
    \label{fig:circuits1}
  \end{figure}
  \begin{enumerate}
  \item[(a)] If the circuit is of dimension $3$ (Cases (A) and (B) in
    Figure \ref{fig:circuits1}), the dual cell is a
    vertex $V$ of $S$ and this vertex is the only singular point.
    \begin{description}
    \item[(a.1)] Either $V$ is adjacent to six edges and nine
      $2$-dimensional polyhedra. Then the dual polytope is IUA-equivalent to a pentatope with
      vertices $(0,0,0)$, $(1,0,0)$, $(0,1,0)$, $(0,0,1)$ and
      $(1,p,q)$ with $p$ and $q$ coprime (Case (A) in Figure~
      \ref{fig:circuits1}, see also Figure \ref{fig:a.1}).
      \begin{figure}[h]
        \centering
        \includegraphics[width=5cm]{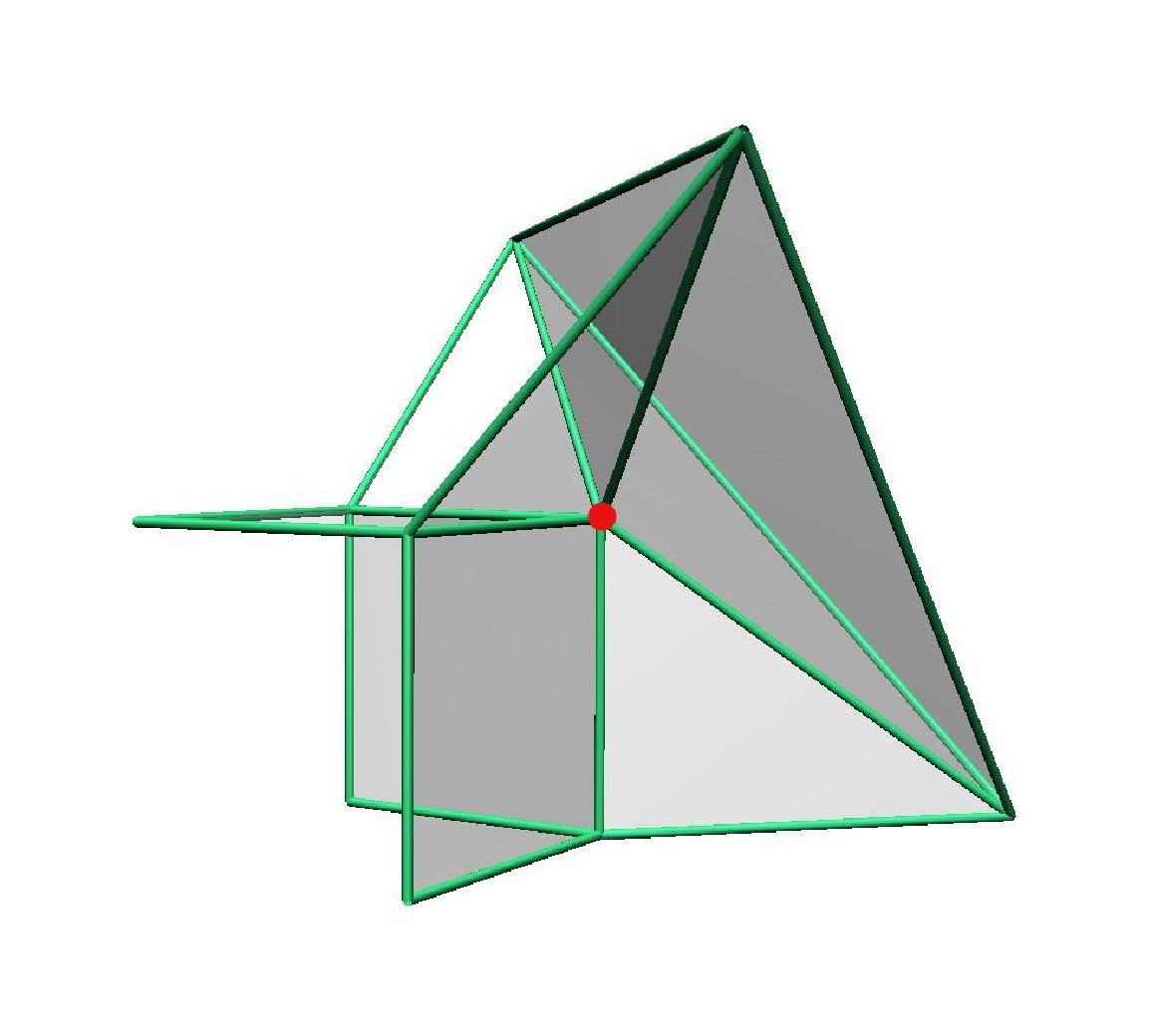}
        \caption{Case (a.1), a singular tropical surface dual to (A) with
          the singular point marked.}
        \label{fig:a.1}
      \end{figure}
    \item[(a.2)] Or $V$ is adjacent to four edges and six
      $2$-dimensional polyhedra, just as a smooth vertex (Case (B) in
      Figure \ref{fig:circuits1}, see also Figure \ref{fig:a.2}). However, if
      we define the multiplicity of a vertex of a tropical
      hypersurface analogously to the case of tropical curves as the
      lattice volume of the corresponding polytope in the dual
      subdivision, then it follows that $V$ is a vertex of higher
      multiplicity. More precisely, the multiplicity can be
      4,5,7,11,13,17,19 or 20. The dual is IUA-equivalent to a tetrahedron with vertices
      $(0,0,0)$, $(1,0,0)$, $(0,1,0)$ and, resp., $(3,3,4)$,
      $(2,2,5)$, $(2,4,7)$, $(2,6,11)$, $(2,7,13)$, $(2,9,17)$,
      $(2,13,19)$, or $(3,7,20)$.
      \begin{figure}[h]
        \centering
        \includegraphics[width=5cm]{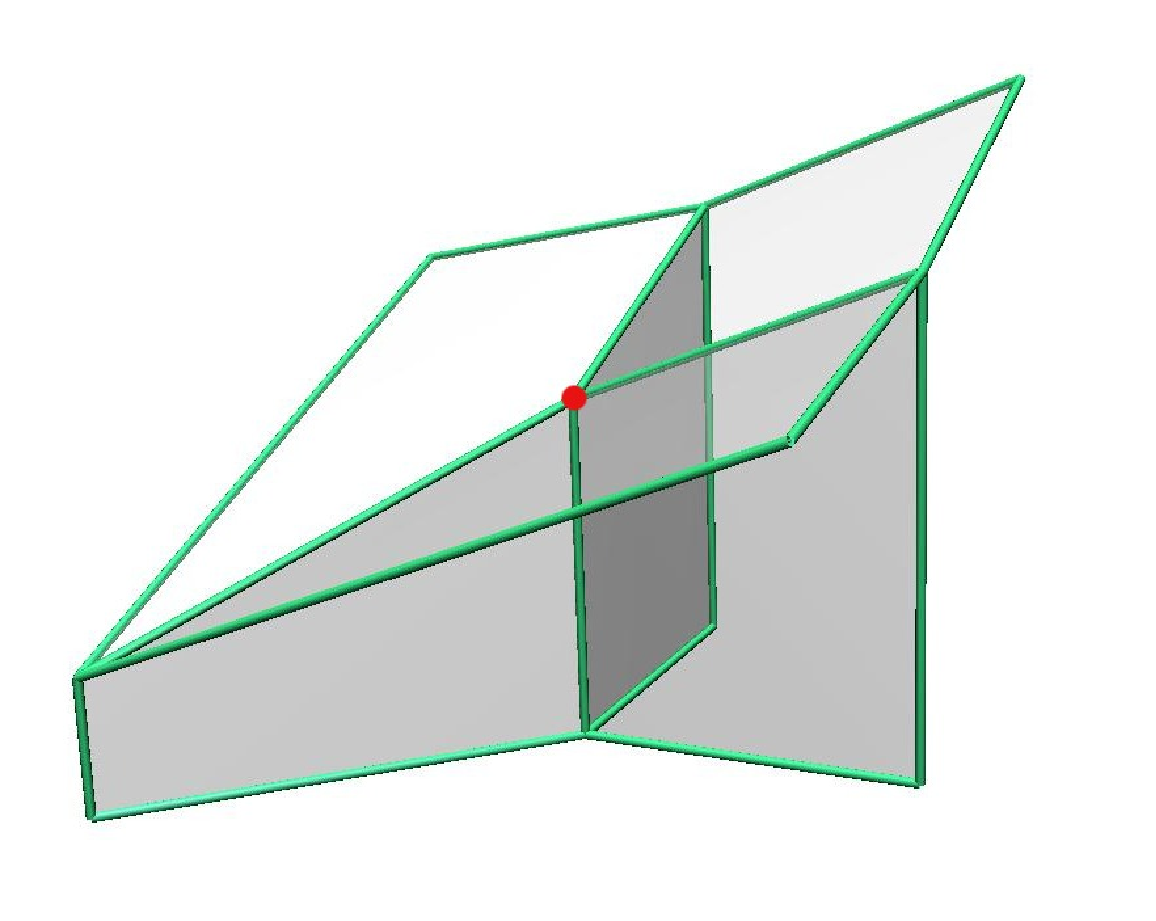}
        \caption{Case (a.2), a singular tropical surface dual to (B) with
          the singular point marked.}
        \label{fig:a.2}
      \end{figure}
    \end{description}
  \item[(b)] If the circuit is of dimension $2$ (Cases (C) and (D) in
    Figure \ref{fig:circuits1}), the dual cell is an
    edge $E$. We have the following cases:
    \begin{description}
    \item[(b.1)] The dual of $E$ is IUA-equivalent to a triangle with vertices $m_a=(0,0,0)$,
      $m_c=(0,1,2)$ and $m_d=(0,2,1)$, i.e.\ $E$ is adjacent to three
      $2$-dimensional cells of $S$ (Case (C) in Figure
      \ref{fig:circuits1}). Each end vertex of $E$ is adjacent
      to four edges and six $2$-dimensional polyhedra, just as a smooth
      vertex.
      \item[\phantom{x}(b.1.1)]
      $E$ is bounded and there is a singularity at the
      midpoint of $E$ or at points which divide $E$ with the ratio
      $3:1$ (see Figure \ref{fig:b.1.1}). Or, $E$ is unbounded and there is a singularity whose distance from the vertex of $E$ depends on six coefficients of the tropical polynomial (see Equation \eqref{eq:b.1.1.4} in Subsection \ref{subsubsec-eq}).
\item[\phantom{x}(b.1.2)]
A bounded edge $E$ admits finitely many (bounded or unbounded) extensions to a \emph{virtual edge} with a singularity at the positions described in (b.1.1) (see Figure \ref{fig:b.1.1.4}). (Subsection \ref{subsubsec-eq} explains the term virtual edge). 
\begin{figure}[h]
        \centering
        \unitlength1cm
        \begin{picture}(14,12)
          \put(0,8.5){\includegraphics[width=5cm]{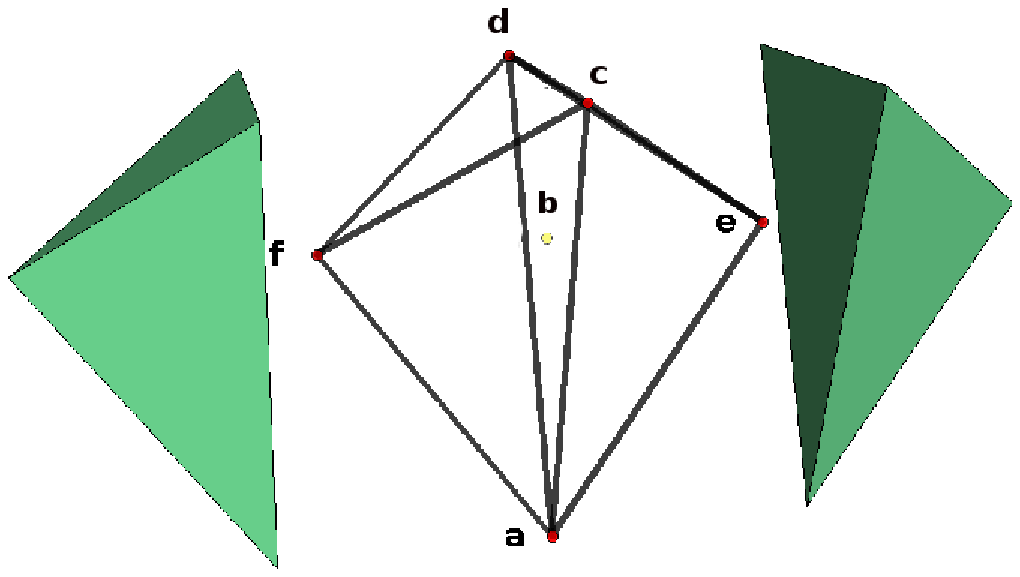}}
          \put(7,8){\includegraphics[width=4cm]{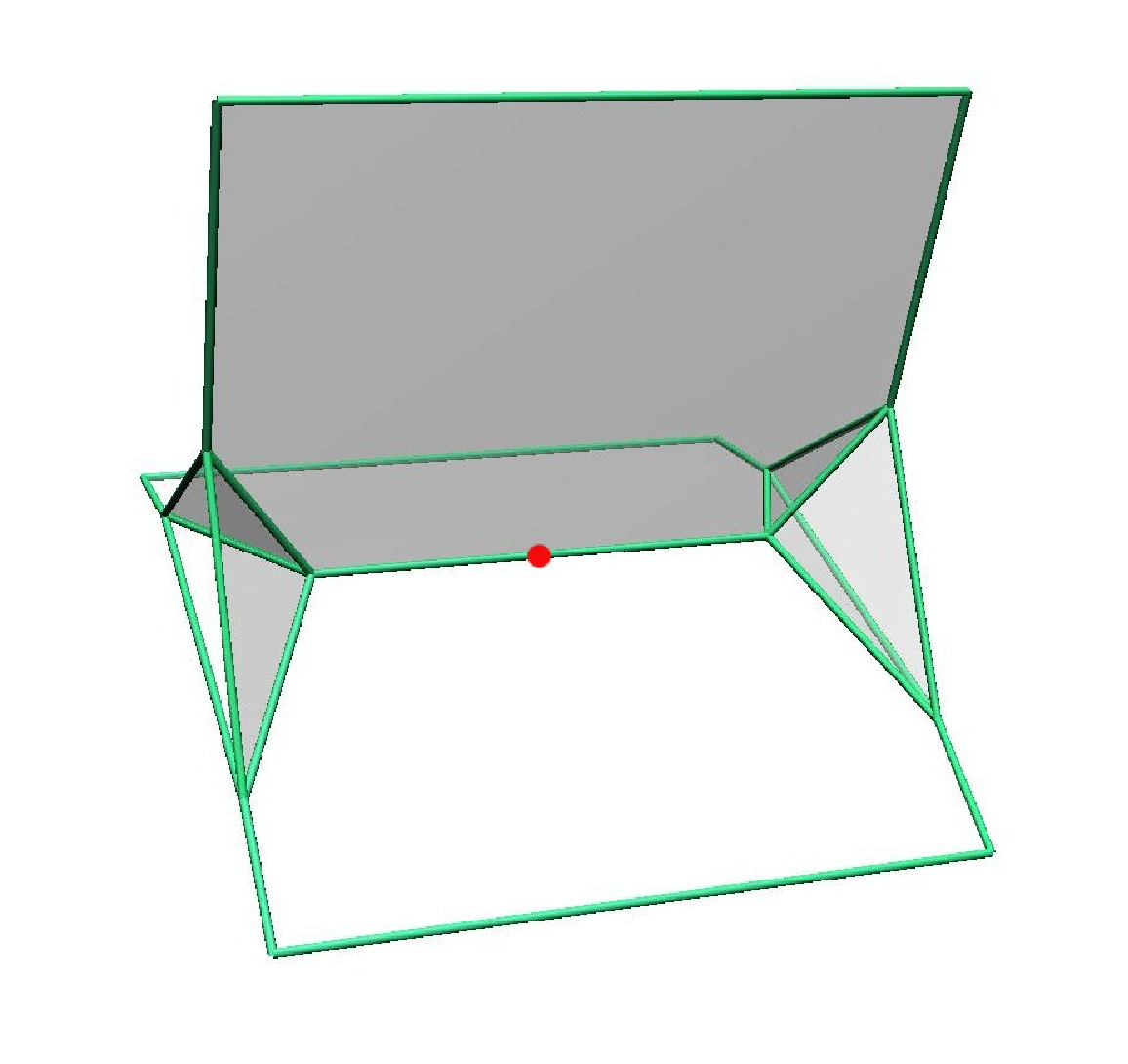}}
          \put(0,0.5){\includegraphics[width=6cm]{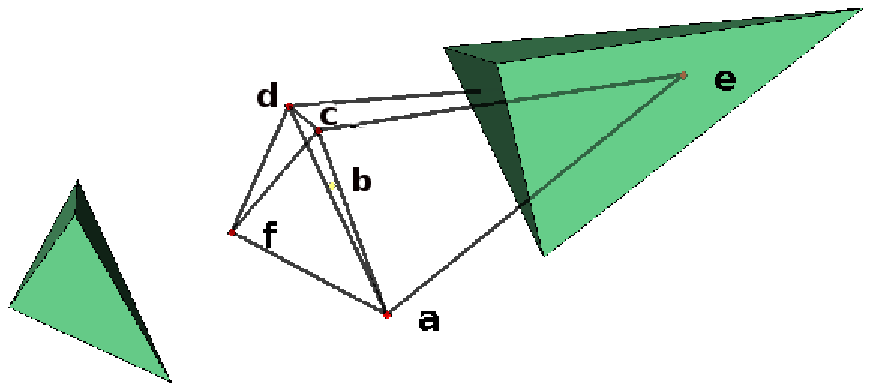}}
          \put(7,0){\includegraphics[width=4cm]{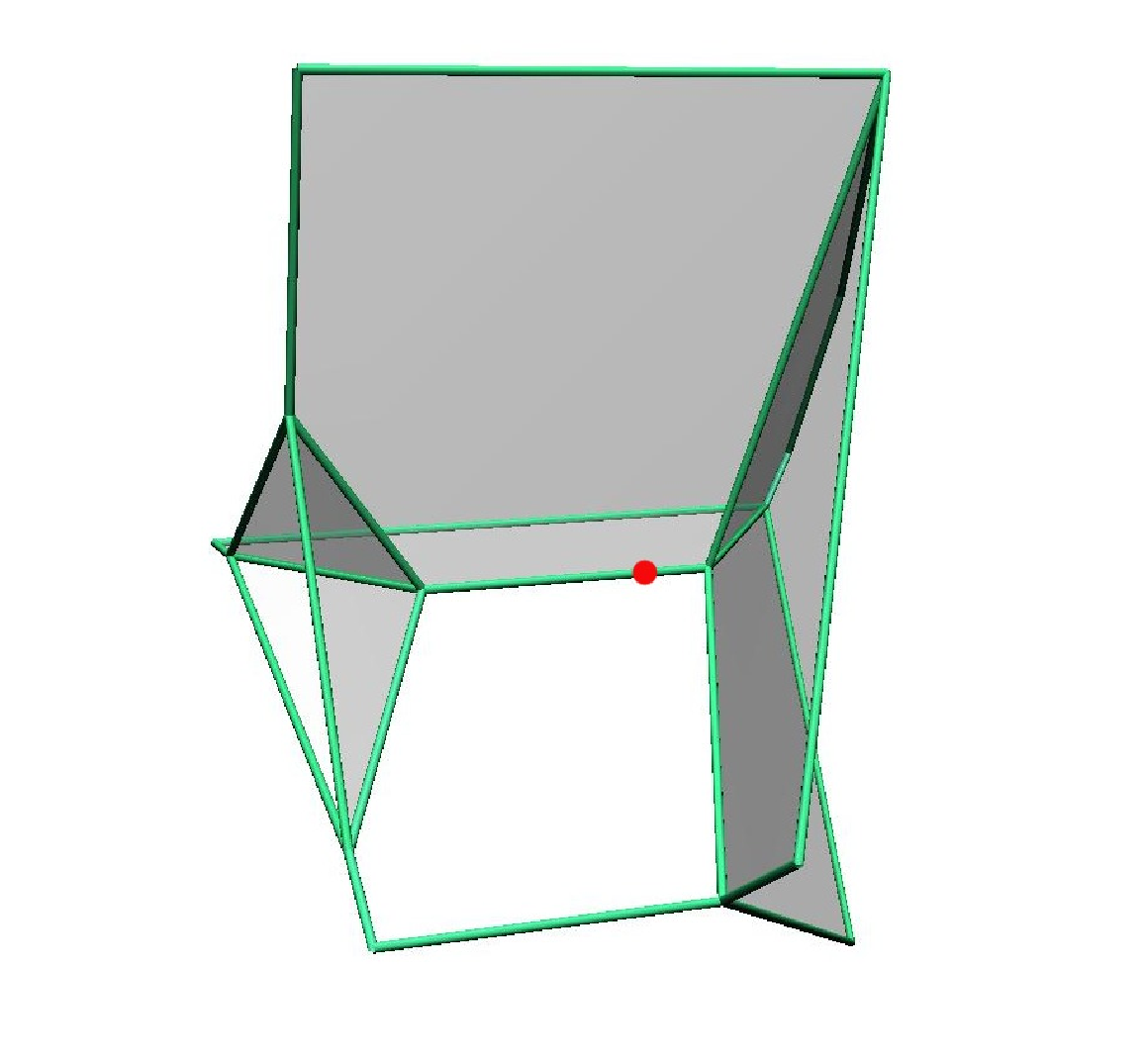}}
        \end{picture}
        \caption{A singular point which divides the edge $E$
          either in the midpoint or with ratio $3:1$, and dual subdivisions (case (b.1.1)).}
        \label{fig:b.1.1}
      \end{figure}

   \begin{figure}[h]
        \centering
        \unitlength1cm
        \begin{picture}(11,5)
          \put(0,1){\includegraphics[width=7.5cm]{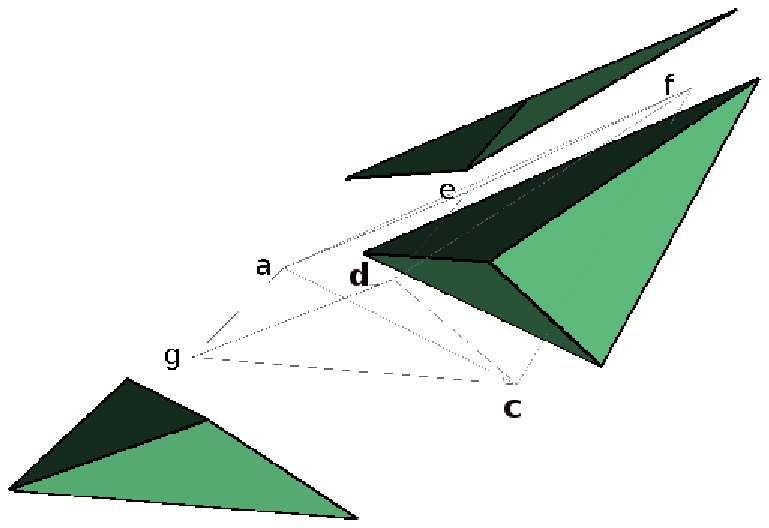}}
          \put(7,0){\includegraphics[width=4.5cm]{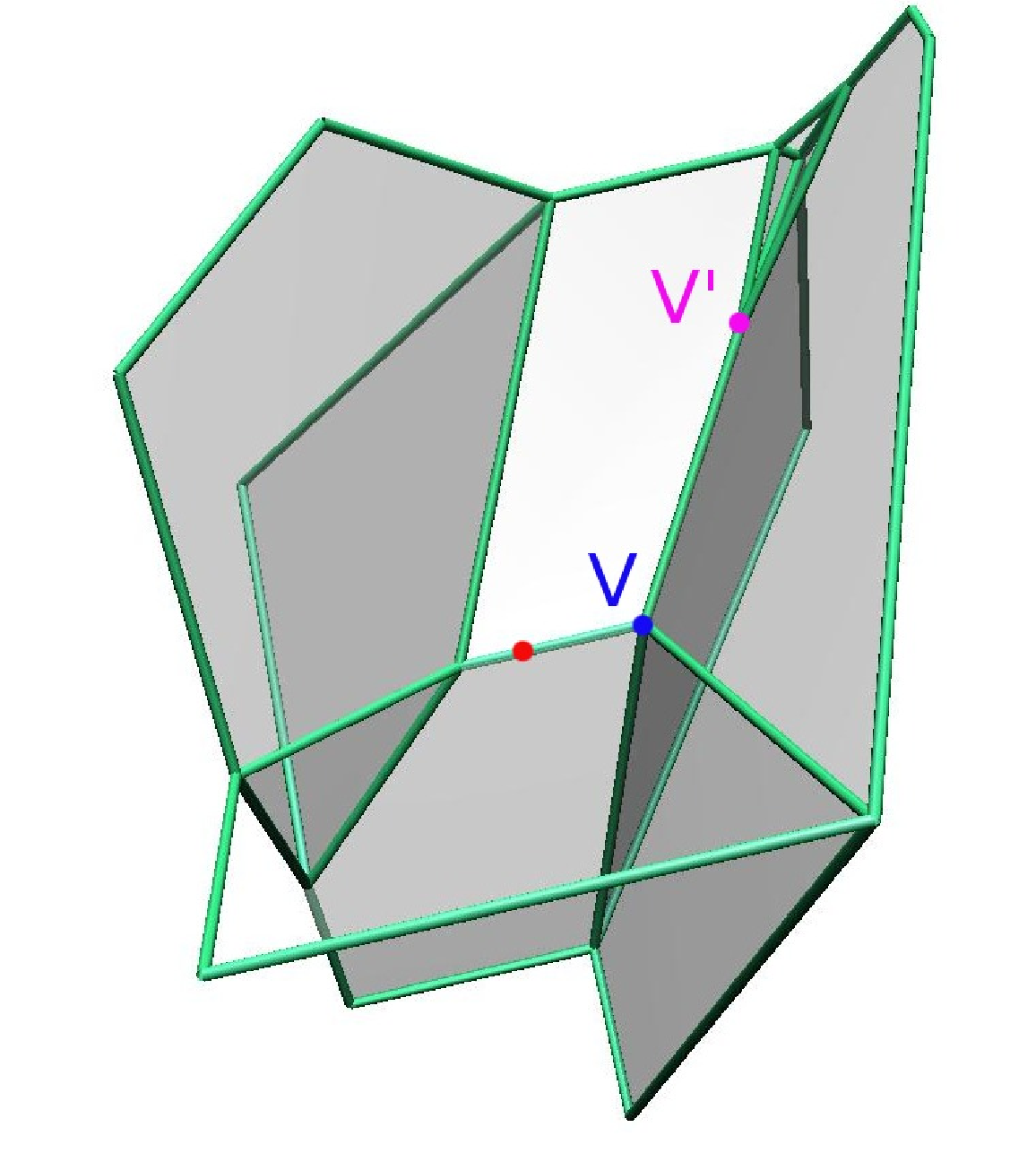}}
        \end{picture}
        \caption{A tropical surface with a singularity as in case (b.1.2) and its dual subdivision. The edge with the singular point can be extended to an unbounded edge containing only the vertex $V$. The singular point is at distance $d$ from $V$, where $d$ depends on the coefficients involving the two vertices $V$ and $V'$.}
        \label{fig:b.1.1.4}
      \end{figure}

    \item[(b.2)] $E$ is dual to a quadrangle, i.e.\ adjacent to four
      $2$-dimensional cells of $S$ (Case (D) in Figure
      \ref{fig:circuits1}, see also Figure \ref{fig:b.2}). $E$ must be bounded and its end
      vertices are each adjacent to five edges and eight $2$-dimensional
      cells. $S$ contains a unique singular point which is the
      midpoint of $E$.
      \begin{figure}[h]
        \centering
        \unitlength1cm
        \begin{picture}(12,4)
          \put(1,1){\includegraphics[width=6cm]{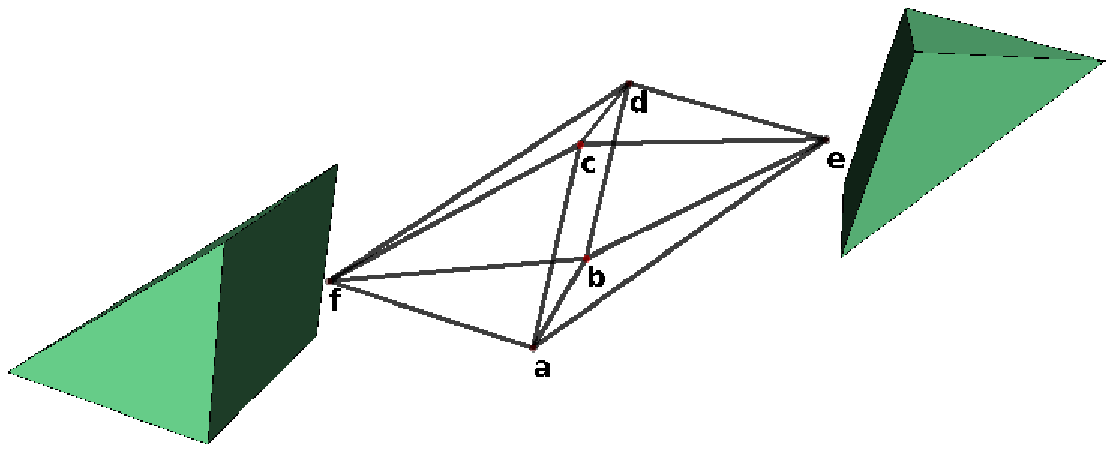}}
          \put(8,0){\includegraphics[width=4.5cm]{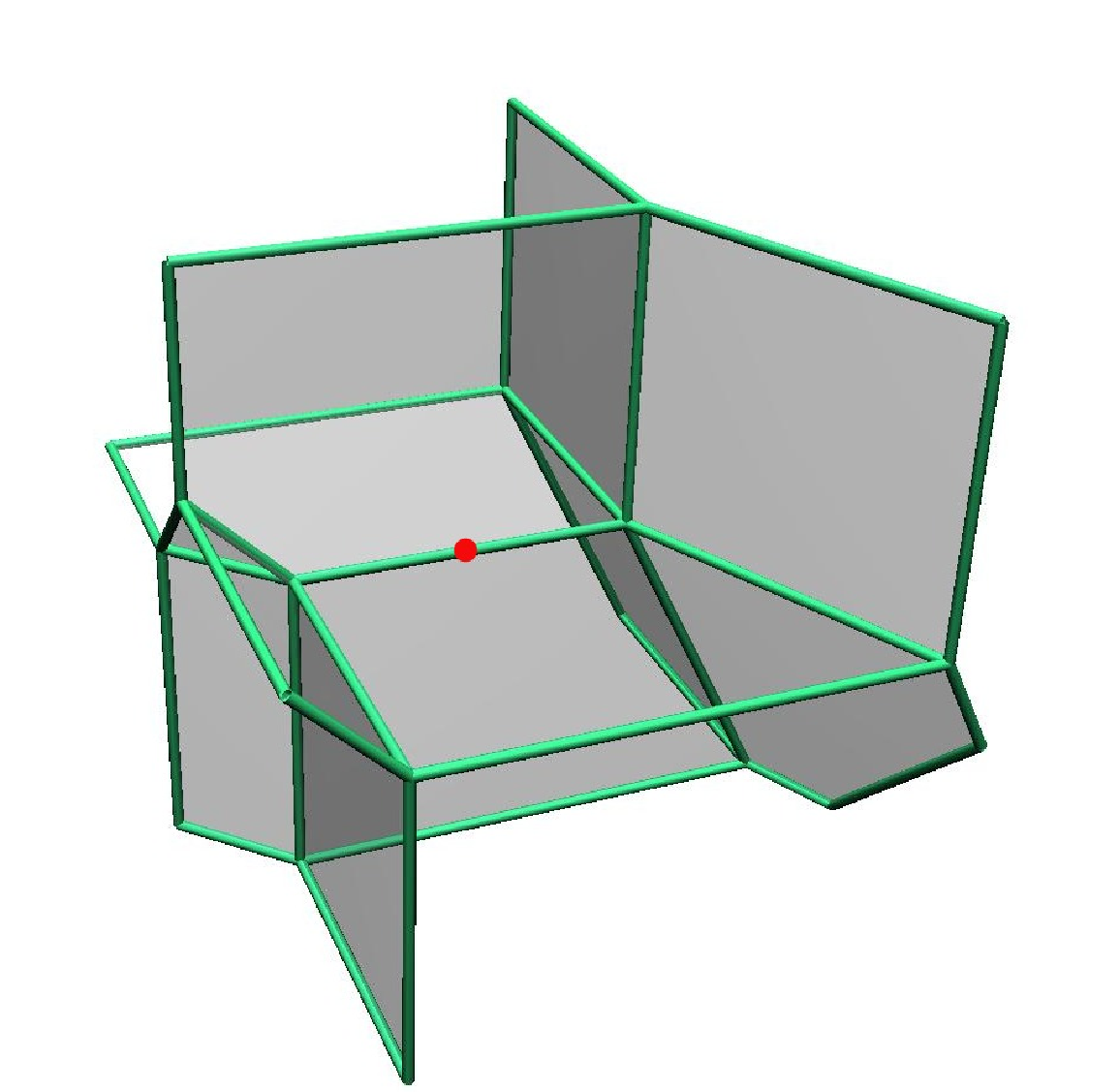}}
        \end{picture}
        \caption{Case (b.2), a dual subdivision with circuit (D) and the
          corresponding singular tropical surface with the singular point marked.}
        \label{fig:b.2}
      \end{figure}
    \end{description}
  \item[(c)] If the circuit is of dimension $1$ (Case (E) in Figure \ref{fig:circuits1}), then the dual is a
    $2$-dimensional cell of $S$. 
 \item[\phantom{x}(c.1)]
If this cell is a triangle or trapeze, there is a singular point at the \emph{weighted barycenter} resp.\ \emph{generalised midpoint}, see Subsection \ref{sec-barycenter} below (for
    an image see Example \ref{ex-thomas}).
 \item[\phantom{x}(c.2)]
An arbitrary $2$-dimensional cell admits finitely many extensions to a triangle or a trapeze, with a singularity at the position described in (c.1).

  \end{enumerate}
\end{theorem}
Subsection \ref{sec-barycenter} referred to in statement (c) of Theorem \ref{thm2}
contains a classification of the possible shapes of the cell dual to
the circuit and explains the terms \emph{weighted barycenter} and
\emph{generalised midpoint}.

\begin{example}\label{ex-thomas}
  A tropical surface $S$ can have
  several singularities, since there may be several singular surfaces
  tropicalising to $S$ with different images for their singular
  point. We give here an example for this behaviour. Consider the
  polynomials
  \begin{displaymath}
    f=(1-3t^5-3t^8)+(-2+t^5)\cdot z +z^2+t^8\cdot \frac{1}{xy}
    +(t^5+t^8)\cdot y +(2t^5+t^8)\cdot x-t^5\cdot x^2yz
  \end{displaymath}
  and
  \begin{displaymath}
    g=(1-3t^6+3t^8)-(2+t^8)\cdot z +z^2+t^8\cdot \frac{1}{xy}
    +(t^5-t^7)\cdot y + (t^5-2t^7)\cdot x +t^5\cdot x^2yz
  \end{displaymath}
  over the field of Puiseux series. They both tropicalise to the
  tropical polynomial
  \begin{displaymath}
    F_u=\max\{0,z,2z,-8-x-y,-5+y,-5+x,-5+2x+y+z\}
  \end{displaymath}
  with $u=(0,0,0,-8,-5,-5,-5)$
  and define thus the same tropical surface $S$. Moreover, both $V(f)$ and
  $V(g)$ are singular, however, $V(f)$ is singular in $(1,1,1)$ which
  tropicalises to $G=(0,0,0)$, while $V(g)$ is singular in $(t,t,1)$
  which tropicalises to $H=(-1,-1,0)$. Thus $S$ has two singular points on
  the quadrangle dual to the circuit formed by $(0,0,0)$, $(0,0,1)$,
  and $(0,0,2)$. The quadrangle is shown in Figure
  \ref{fig:twobarycenters}, and $G=\frac{1}{3}\cdot(A+B-E)$ and
  $H=\frac{1}{3}\cdot(C+D+E)$ are weighted barycentes of the vertices
  $A$ and $B$ respectively $C$ and $D$ with the virtual vertex $E$ in
  the sense of Theorem \ref{thm2} and Subsection \ref{sec-barycenter}
  (see also Remark \ref{rem:ex-thomas} and Examples
  \ref{ex:barycenter1} and \ref{ex:barycenter2}).
  \begin{figure}[h]
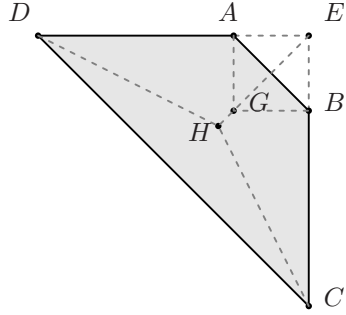

    \centering
    \begin{texdraw}
      \drawdim cm \relunitscale 0.2
      \move (5 -13) \lvec (-13 5) \lvec (0 5) \lvec (5 0)
      \lvec (5 -13) \lfill f:0.9
      \lpatt (0.3 0.5)
      \setgray 0.5
      \move (5 0) \lvec (5 5) \lvec (0 5)
      \move (5 0) \fcir f:0 r:0.2
      \move (0 5) \fcir f:0 r:0.2
      \move (5 5) \fcir f:0 r:0.2
      \move (-13 5) \fcir f:0 r:0.2
      \move (5 -13) \fcir f:0 r:0.2
      \move (0 0) \fcir f:0 r:0.2
      \move (-1 -1) \fcir f:0 r:0.2
      \htext (-15 6) {$D$}
      \htext (6 -13) {$C$}
      \htext (-1 6) {$A$}
      \htext (6 0) {$B$}
      \htext (6 6) {$E$}
      \htext (1 0) {$G$}
      \htext (-3 -2) {$H$}
      \move (-13 5) \lvec (-1 -1) \lvec (5 -13)
      \move (0 5) \lvec (0 0) \lvec (5 -0)
      \move (5 5) \lvec (-1 -1)
    \end{texdraw}
    \caption{Two singular points on a tropical surface as weighted barycenters.}
    \label{fig:twobarycenters}
  \end{figure}
\end{example}

Theorem \ref{thm2} gives necessary conditions for the geometry of a
singular tropical surface. We can also formulate a sufficient
condition, which follows from Lemma \ref{lem-chains} and the classification:
\begin{theorem}
 \label{thm3}
 Let $F_u=\max_{(i,j,k)\in\mathcal{A}}\{u_{(i,j,k)}+ix+jy+kz\}$
  define a tropical surface $S$ dual to a
  marked subdivision  of maximal-dimensional geometric type. Assume that the dual subdivision corresponds to
  a cone of codimension $1$ in the secondary fan, and its unique circuit is one of the IUA-types shown in Figure \ref{fig:circuits1}.
Then $S$ is a singular tropical surface if
\begin{itemize}
\item either the circuit is of type (A) or (B),
\item or the circuit is of type (C), (D), or (E) and it does not lie on the boundary of $\Delta$,
\item or the circuit is of type (C), lies on $\partial\Delta$
and it is the base of a pyramid $P$ of the dual subdivision such that $\mbox{vol}(P)=9$ and $P\subsetneqq\Delta$,
\item or the circuit is of type (E), lies on $\partial\Delta$, and $\Delta$ contains 3 more points as described in Proposition
\ref{prop:classification(E)1} below, or 4 more points as described in Section \ref{subsec:circuitE} below. \end{itemize} 

Furthermore, let $p\in S$ be a point in the cell dual to the circuit, and assume
  $p$ satisfies conditions (a), (b) or (c) of Theorem \ref{thm2}
  above. Then $S$ is the tropicalization of an algebraic surface with a
  singularity tropicalizing to $p$ if and only if after shifting $S$
  such that $p$ becomes the origin (and accordingly adding lineality
  vectors to the coefficients $u$ such that they become equal along the
  circuit, see Section \ref{sec:tropicaldiscriminant}) the flag of
 subsets $\mathcal{F}(u)$ (see Subsection \ref{subsec-family})
  either is a flag
  satisfying the conditions of Lemma \ref{lem-chains}, or is in the
  boundary of such a flag.
\end{theorem}

Note that Theorems \ref{thm2} and \ref{thm3} together give a complete classification of maximal-dimensional geometric tropical surfaces and their singular points, and both the necessary and sufficient criteria are easy to verify in any concrete example. 

 For circuits of type (C), the singularity condition may impose non-local geometric conditions. Non-local here means that they involve cells of the tropical surface which are not faces of the cell dual to the circuit.
The following example presents such a situation.

\begin{example}\label{ex-threepoints}
  Let us consider the point configuration $\mathcal{A}$ with $m_a=(0,0,0)$,
  $m_b=(0,1,1)$, $m_c=(0,1,2)$, $m_d=(0,2,1)$, $m_e=(1,1,1)$,
  $m_f=(3,0,2)$ and $m_g=(-1,1,0)$, and a tropical surface $S$ defined
  by a tropical polynomial $F_u$, $u=(u_a,u_b,\ldots,u_g)$.
We assume that $u_a=u_b=u_c=u_d\geq u_e,u_f,u_g$, or equivalently,
we assume that the edge $E$ dual to the circuit satisfies $y=z=0$.
From Theorem \ref{thm2} and Subsection \ref{subsubsec-eq} can conclude that
$S$ can be singular at the point $p$ which divides $E$ with ratio
$3:1$, or at the point $q$ whose position is determined by Equation
  \eqref{eq:b.1.1.4}, or at the point $r$ with coordinates $(\frac{u_g-u_e}{2},0,0)$. (The point $r$ is the midpoint of an extension of $E$ to a virtual edge, see Subsection \ref{subsubsec-eq}.)
 For the points $q$ and $r$, the position of the
  singular point is not \emph{locally} determined, i.e.\ it is not determined
  purely by the linear forms in $F_u$ corresponding to the part of the
  subdivision which is dual to the edge $E$ and its end points, but
  it involves the vertex $V'$ of $S$ determined by the
  polytope $m_a,m_d,m_e,m_f$ (see also Figure \ref{fig:b.1.1.4}).

We now want to specify the sufficient conditions we observe in
Theorem \ref{thm3} in this situation in order to decide which of the
points $p$, $q$ or $r$ is a singular point of the tropical surface. 
If we
move $p$ to the origin, this corresponds to adding the vector
$\frac{u_e-u_f}{2}\cdot (0,0,0,0,1,3,-1)$ to the coefficient vector
$(u_a,u_b,u_c,u_d,u_e,u_f,u_g)$. The new coefficients satisfy the
conditions of Lemma \ref{lem-chains} if and only if the new
$g$-coefficient is smaller than the new $e$ and $f$-coefficients
which became equal. This is the case if and only if $2u_e> u_g+u_f$.
Moving the point $r$ to the origin corresponds to adding the vector $\frac{u_g-u_e}{2}\cdot (0,0,0,0,1,3,-1)$. The new coefficients satisfy the conditions if and only if the new $f$-coefficient which equals is smaller than the new $g$ and $e$-coefficients, which again is the case if and only if $2u_e> u_g+u_f$.
Thus $S$ is singular at both points $p$ and $r$ if and only if
$2u_e\geq u_g+u_f$.

Moving the point $q$ to the origin corresponds
to adding the vector $\frac{u_g-u_f}{2}\cdot (0,0,0,0,1,3,-1)$ to
the coefficient vector. The new coefficient vector satisfies the
conditions of Lemma \ref{lem-chains} if and only if $2u_e< u_g+u_f$.
Thus $q$ is a singular point of $S$ if and only if $2u_e\leq
u_g+u_f$. If $2u_e=u_g+u_f$ then $q=r$ and the coefficient vector is
in the boundary of three weight classes satisfying the conditions of
Lemma \ref{lem-chains}. In any case, we either have one or two singular points, depending on the coefficients of $u$.

\end{example}

\begin{remark}\label{rem-discrim}
The classification is closely related to the study of
$\Delta$-equivalence classes of marked subdivisions (see Section 11.3
of \cite{GKZ}), since by Theorem 1.1 of \cite{DFS05}, the tropical
discriminant (which equals the codimension one subfan of the secondary fan that
groups maximal dimensional cones of the secondary fan into
$\Delta$-equivalence classes) equals the Minkowski sum of the tropicalisation of the family
of curves with a singularity in a fixed point and its lineality space.
This explains why the dual marked subdivisions of maximal-dimensional geometric singular
tropical surfaces correspond to codimension one cones of the
secondary fan which  separate two non-$\Delta$-equivalent maximal
cones (see 11.3.10 of \cite{GKZ} for the smooth case):
understanding the combinatorial types of singular tropical
hypersurfaces is equivalent to understanding $\Delta$-equivalence
classes.
Since understanding $\Delta$-equivalence classes combinatorially is an open
problem for dimension larger than $3$, this connection restricts further generalizations of Theorems \ref{thm2} and \ref{thm3} to
higher dimensions.
\end{remark}


This paper is organised as follows. In Section \ref{sec-basic} the
basic notions will be introduced, most prominently the
tropicalisation $\Trop(\Ker(A))$ of the family of surfaces in a
given toric threefold which are singular at $(1,1,1)$. We also
explain how $\Trop(\Ker(A))$ comes in a natural way with a fan
structure induced by the matroid associated to $A$, and we describe
the full-dimensional cones of this fan as weight classes associated
to flags of flats  (see Lemma \ref{lem-chains}). It is well known
from \cite{DFS05} that the secondary fan of the point configuration
corresponding to $A$  is the Minkowski sum of $\Trop(\Ker(A))$ and
the lineality space. In Section \ref{sec:tropicaldiscriminant} we
reconsider how the Minkowski sum of a cone in $\Trop(\Ker(A))$ with
the lineality space can lie in cones of the secondary fan, and we
use this to introduce the notion of a generic singular surface as
well as to prove Theorem \ref{thm1}. Section \ref{sec-class} is
devoted to the classification of generic singular tropical surfaces
of maximal-dimensional geometric type, and the classification works along the
classification of weight classes in Lemma \ref{lem-chains}. For the
classification also polytopes with certain properties have to be
classified, and the corresponding classification results can be
found in Section \ref{sec-class} too.

\subsection{Acknowledgements}
We would like to thank Christian Haase for useful discussions. The
images were obtained with the aid of Polymake \cite{GJ11},
Javaview \cite{javaview}, jReality \cite{jreality}, tropicalinsect
\cite{tropicalinsect}, xfig and texdraw.
The authors were supported by the
  Hermann-Minkowski Minerva Center for Geometry at the Tel Aviv
  University, and by the DFG-grant MA 4797/3-1 as part of the priority
  program SPP 1489. The third author was also supported by the Israeli
  Science Foundation grant no. 448/09. We would like to thank an anonymous referee for valuable comments on a first draft of this paper.

\section{Notations and basic facts}\label{sec-basic}
In this section, we fix notations and collect basic properties of
the family of surfaces with a singularity in a fixed point and its
tropicalisation, the Bergman fan of the corresponding linear ideal.
The content of this section is parallel to Sections 1, 2 and 3.1 of
\cite{MMS09},
only now we deal with surfaces instead of curves. We
omit proofs in this section, since they are all straightforward
generalisations of the corresponding statements in \cite{MMS09}.

\subsection{The family of surfaces with a singularity in a fixed point}\label{subsec-sing}
   Fix a non-degenerate convex
   lattice polytope  $\Delta\subset\R^3$ and denote by
   $\mathcal{A}=\Delta\cap\Z^3=\{m_1,\ldots,m_s\}$ the lattice points of $\Delta$. For
   any field $\K$ there is a toric threefold
   $\Tor_\K(\Delta)$ associated to $\Delta$ and it comes with
   the tautological line bundle ${\mathcal L}_\Delta$
   generated by the global sections $\{x^iy^jz^k\ :\
   (i,j,k)\in\mathcal{A}\}$. The torus $(\K^*)^3$ is embedded in
   $\Tor_\K(\Delta)$ via
   \begin{displaymath}
     \Psi_{\mathcal{A}}:(\K^*)^3\longrightarrow \TP_\K^{\mathcal{A}}:(x,y,z)\mapsto
     \big(x^iy^jz^k\;|\;(i,j,k)\in\mathcal{A}\big)
   \end{displaymath}
   and inside the torus the elements in the linear system $|{\mathcal
     L}_\Delta|$ are defined by the equations
   \begin{displaymath}
     f_{a}=\sum_{(i,j,k)\in\mathcal{A}}a_{(i,j,k)}\cdot x^i\cdot
     y^j \cdot z^k=0
   \end{displaymath}
   with $a=(a_{(i,j,k)}\;|\;(i,j,k)\in\mathcal{A})\in(\TP_K^{\mathcal A})^*$.
   $|{\mathcal L}_\Delta|$
   contains a nonempty linear subsystem $\Sing_\bp(\Delta)$ of surfaces with
   a singularity at the point $\bp=(1,1,1)$. The equations for this
   subsystem are the linear equations
   \begin{displaymath}
     f_{a}(\bp)=0,\;\; \frac{\partial f_{a}}{\partial
       x}(\bp)=0,\;\;\frac{\partial f_{a}}{\partial
       y}(\bp)=0,\;\;\frac{\partial f_{a}}{\partial
       z}(\bp)=0,
   \end{displaymath}
   or equivalently we can say that the family $\Sing_\bp(\Delta)$ is the
   kernel of the $4\times s$ matrix
   \begin{displaymath}
     A=\left(\begin{matrix}
         1&\ldots &1\\
         m_1&\ldots&m_s
       \end{matrix}\right).
   \end{displaymath}
   Notice that $A$ is just the matrix of the point configuration $\mathcal{A}$,
   after raising the points to the $\{t=1\}$-plane in $\R^4$, if we
   choose the coordinates $(t,x,y,z)$ on $\R^4$.

\subsection{Tropicalisations}\label{subsec-trop}
Let $\K$ denote the field of Puiseux series and $\val$ the valuation sending a Puiseux series to the smallest exponent.
\longer{For tropicalisations, we use an algebraically
   closed field $\K$ with a non-archimedean valuation
   $\val:\K^*\longrightarrow\R  $
   satisfying $\val(\Z)=0$ and whose value group is dense in $\R$, e.g.\
   the algebraic closure $\overline{\C(t)}$ of the field of rational
   functions over $\C$, or $\C\{\{t\}\}$ the field of Puiseux series,
   or a field of generalised Puiseux series as in \cite{Mar08}. In
   each of these cases the elements of the field can be represented
   by generalised power series of the form
   \begin{displaymath}
     p=a_1 t^{q_1}+a_2 t^{q_2}+\ldots
   \end{displaymath}
   with complex coefficients and real exponents, and the valuation
   maps $p$ to the least exponent $q_1$ whose coefficient $a_1$ is
   non-zero.}
 For an ideal $I\subset \K[x_1^\pm,\ldots,x_n^\pm]=\K[\bx^\pm]$ determining a variety
   $V=V(I)\subset (\K^\ast)^{n}$ we define the \emph{tropicalisation} of $V$ to be
   \begin{displaymath}
     \Trop(V):= \overline{\{(-\val(x_1),\ldots,-\val(x_n)) \;|\;
       (x_1,\ldots,x_n)\in V(I) \}},
   \end{displaymath}
   i.e.\ we map $V$ componentwise with the negative of the valuation
   map and take the topological closure in $\R^n$.
  \longer{ If the ideal $I$ is generated by homogeneous polynomials we
   may alternatively consider $V(I)$ inside $\TP_\K^n$ and we, consequently,
   should consider $\Trop(V)$ modulo the linear space
   spanned by $(1,\ldots,1)$, i.e.\ we
   should identify $\Trop(V)$ with its image in $\R^n / (1,\ldots,1)$.}

   We consider tropicalisations in two  situations:
   \begin{itemize}
   \item \emph{The tropicalisation of $\Sing_\bp(\Delta)=\ker(A)$:} The
     linear space $V=\ker(A)$ is defined by linear equations over $\Q$. $\Trop(V)$  is the so-called  \emph{Bergman fan} of $I$ (\cite{FS05},
     \cite{AK06}). We will study the Bergman fan $\Trop(\ker(A))$
     further in Subsection \ref{subsec-family}. Note, since the linear
     generators of $A$ are homogeneous, we will consider $\Trop(V)$
     modulo the vector space spanned by $(1,\ldots,1)$.
     That is, we consider $\Trop(\Sing_{\bp}(\Delta))=\Trop(\Ker(A))$
     as a fan in $\R^{s-1}=\R^{\mathcal{A}}/(1,\ldots,1)$.
   \item \emph{The tropicalisation of a surface $V(f_{a})$ with $a\in
       \Sing_\bp(\Delta)$:} This is an example of a tropical
     hypersurface. If $V$ is a hypersurface defined by $f=\sum a_{m}
       {\bx}^{m}$, then its tropicalisation equals the locus of
     non-differentiability of the \emph{tropical polynomial}
     \begin{displaymath}
       \trop f:\R^n\longrightarrow\R:{\bx}\mapsto \max \{ -\val(a_{m})+m\cdot {\bx}\}
     \end{displaymath}
     by Kapranov's Theorem (see
     \cite[Theorem~2.1.1]{EKL06}).
   \end{itemize}
Let us first study the hypersurface case more closely.

\subsection{Tropical hypersurfaces and dual marked subdivisions}\label{subsec-hypersurf}
Tropical hypersurfaces are dual to \emph{marked subdivisions}
$T=\{(Q_1,\mathcal{A}_1),\ldots,(Q_k,\mathcal{A}_k)\}$ (where the
$Q_i$ are
polytopes and the $\mathcal{A}_i$ marked integer points,
see \cite[Definition 7.2.1]{GKZ} resp.\ \cite[Section 2]{MMS09}).
      We define the \emph{type} of a marked subdivision to be the
     subdivision, i.e.\ the collection of $Q_i$ without the
     markings.

   For a finite subset $\mathcal{A}$ of the lattice $\Z^d$ we denote
   by $\R^{\mathcal{A}}$ the set of vectors indexed by the lattice
   points in $\mathcal{A}$.
   A point $u\in\R^{\mathcal{A}}$ induces a \emph{regular (or coherent)} marked subdivision
   of $\Delta$ by considering the convex hull of
   \begin{equation}\label{e1}
     \big\{(m,u_m)\;\big|\;m\in\mathcal{A}\}\subset\R^d\times \R
   \end{equation}
   in $\R^{d+1}$, and projecting the upper faces onto $\R^d$. A
   lattice point $m$ is marked if the point $(m,u_{m})$ is
   contained in one of the upper faces. We say two points $u$ and $u'$ in $\R^{\mathcal{A}}$
   are equivalent if and only if they induce the same regular marked subdivision
   of $\Delta$.  This defines an equivalence relation on
   $\R^{\mathcal{A}}$ whose equivalence
   classes are the relative interiors of convex cones. The collection of these cones is
   the \emph{secondary fan} of $\Delta$.

   Regular marked subdivisions of $\Delta$ are dual to tropical hypersurfaces
   (see e.g.\ \cite[Prop.\ 3.11]{Mi03}).
   Given a point $u\in\R^{\mathcal{A}}$
   it defines a tropical hypersurface $S_F$ as the locus of
   non-differentiability of the tropical polynomial \label{page:troppoly}
   \begin{displaymath}
     F_u=\max\{u_{m}+m\cdot {\bx}\;|\; m \in\mathcal{A}\},
   \end{displaymath}
   and it defines a regular subdivision of $\Delta$.
Each $k$-dimensional polytope in the subdivision is dual to a
$d-k$-dimensional orthogonal polyhedron of the tropical hypersurface.

For tropical surfaces dual to a marked subdivision of a polytope in $\R^3$, this means more precisely:
\begin{itemize}
 \item each $3$-dimensional polytope in the subdivision is dual to a vertex of the tropical surface;
\item each $2$-dimensional face in the subdivision is dual to an edge of the tropical surface, which
is perpendicular to the plane spanned by the $2$-dimensional face;
\item each edge of the subdivision is dual to a perpendicular $2$-dimensional polyhedron of the tropical surface. The weight of a
$2$-dimensional polyhedron of the tropical surface is defined to be $\#(e\cap \Z^3)-1$, where $e$ is the dual edge in the marked subdivision.
\end{itemize}

   The duality implies that we can deduce the type of the marked
   subdivision from the tropical hypersurface $S_F$, but not the
   markings. To deduce the markings, we need to know the
   coefficients $u_{m}$.

   Obviously, the vector $(1,\ldots,1)$ is contained in the lineality space of the secondary fan.
   Therefore we can mod out this vector and consider the resulting fan
   in $\R^{s-1}= \R^{\mathcal{A}}/(1,\ldots,1)$ with $s=\#\mathcal{A}$.
   We have seen above that every point $u$ in $\R^{\mathcal{A}}$
   defines a tropical hypersurface via the tropical polynomial
   $F_u=\max\{u_{m}+m\cdot {\bx}\}$. Of course, adding $1$ to each
   coefficient $u_{m}$ does not change the tropical hypersurface associated
   to this polynomial. Hence if we consider $\R^{\mathcal{A}}$ as
   a parametrising space for tropical hypersurfaces, it makes sense to mod
   out the linear space spanned by $(1,\ldots,1)$, and we will do so in what follows. By abuse
   of notation, we call the fan in $\R^{s-1}$ that we get from the
   secondary fan in this way also the \emph{secondary fan}.

     The identification of $\R^{\mathcal{A}}$ with $\R^s$,
     $s=\#\mathcal{A}$, is done by fixing an ordering of the elements
     of $\mathcal{A}$, say $m_1,\ldots,m_s$. When referring to an
     element $u\in \R^{\mathcal{A}}=\R^s$ we will sometimes refer to
     the coordinates of $u$ as $u_m$ with $m\in\mathcal{A}$ and
     sometimes simply as $u_i$ with $i=1,\ldots,s$. This should not
     lead to any ambiguity.

  \longer{ Let $T=\{(Q_l,\mathcal{A}_l)\;|\;l=1,\ldots,k\}$, be a marked
   subdivision of $\Delta$.
   \begin{displaymath}
     \ker(A)=\Big\{(a_{m})\in \R^{\mathcal{A}} \;\Bigm|\; \sum\nolimits_{m}
     a_{m}\cdot m=0, \;\sum\nolimits_{m} a_{m}=0\Big\}
   \end{displaymath}
   is the space of affine relations among the lattice points $m$ of $\Delta$.
   For any $l$, let
   \begin{displaymath}
     L_{\mathcal{A}_l}= \{(a_{m})\in L\;|\;a_{m}=0
     \mbox{ for } m\notin \mathcal{A}_l\}
   \end{displaymath}
   be the space of affine relations among the elements of
   $\mathcal{A}_l$, and let $L_T$ be their sum.

   \begin{lemma}
     The codimension of the cone of the secondary fan corresponding to
     the marked subdivision $T$ equals $\dim(L_T)$.

     In particular, a cone in the secondary fan corresponding to a
     marked subdivision is top-dimensional if and only if the marked
     subdivision is a triangulation, i.e. all polytopes $Q_i$ are
     triangles and in each $Q_i$ no other point besides the vertices
     is marked.
   \end{lemma}
For a proof, see \cite[Corollary 2.7]{GKZ}.}

     \begin{remark}\label{rem-circuit}
     A cone in the secondary fan of codimension one contains
     exactly one circuit, i.e.\ a set of
     lattice points that is affinely dependent but such that each
     proper subset is affinely independent.
     A circuit in $3$-space consists either of the $5$ vertices of a
     pentatope such that each subset of $4$ vertices spans the space (A), or of
     the four vertices of a simplex and an interior point (B), or of $4$ points
     in a plane as in (C) resp.\ (D), or of $3$ points on a line (E), as depicted in Figure
     \ref{fig:circuits}.
     \begin{figure}[h]
       \centering
\begin{picture}(0,0)%
\includegraphics{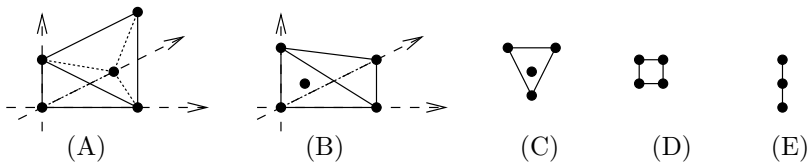}%
\end{picture}%
\setlength{\unitlength}{3947sp}%
\begingroup\makeatletter\ifx\SetFigFont\undefined%
\gdef\SetFigFont#1#2#3#4#5{%
  \reset@font\fontsize{#1}{#2pt}%
  \fontfamily{#3}\fontseries{#4}\fontshape{#5}%
  \selectfont}%
\fi\endgroup%
\begin{picture}(5042,992)(3064,-2966)
\put(7876,-2911){\makebox(0,0)[lb]{\smash{{\SetFigFont{10}{12.0}{\familydefault}{\mddefault}{\updefault}{\color[rgb]{0,0,0}(E)}%
}}}}
\put(3451,-2911){\makebox(0,0)[lb]{\smash{{\SetFigFont{10}{12.0}{\familydefault}{\mddefault}{\updefault}{\color[rgb]{0,0,0}(A)}%
}}}}
\put(4951,-2911){\makebox(0,0)[lb]{\smash{{\SetFigFont{10}{12.0}{\familydefault}{\mddefault}{\updefault}{\color[rgb]{0,0,0}(B)}%
}}}}
\put(6301,-2911){\makebox(0,0)[lb]{\smash{{\SetFigFont{10}{12.0}{\familydefault}{\mddefault}{\updefault}{\color[rgb]{0,0,0}(C)}%
}}}}
\put(7126,-2911){\makebox(0,0)[lb]{\smash{{\SetFigFont{10}{12.0}{\familydefault}{\mddefault}{\updefault}{\color[rgb]{0,0,0}(D)}%
}}}}
\end{picture}%
       \caption{Circuits in $3$-space.}
       \label{fig:circuits}
     \end{figure}
   \end{remark}

\begin{definition}\label{def-mdg-type} Given a tropical surface $S$, we have seen above that it is dual to a type $\alpha=\{Q_1,\ldots,Q_k\}$ of a marked subdivision. We call $\alpha$ also the \emph{type of the tropical surface}.
We can parametrise all tropical surfaces of a given type by an
unbounded polyhedron in $\R^{3\cdot v}$, where $v$ denotes the number of
vertices of $S$. We associate a point in $\R^{3\cdot v}$ to a
tropical surface by collecting all coordinates of vertices.
The polyhedron is defined by equations and inequalities that we can
deduce from the type and that tell us which vertices are connected by
an edge of which direction. We define the
\emph{dimension} $\dim(\alpha)$ of a type $\alpha$ to be the
dimension of this parametrising polyhedron.
If the tropical surface $S$ is singular and
$\dim(\alpha)=\#(\Delta\cap\Z^3)-2$, which is the maximal
possible value for singular tropical surfaces with Newton polytope
$\Delta$, we say that $S$ is of \emph{maximal-dimensional geometric type}.
\end{definition}

   For the following lemma recall that we consider the secondary fan
   of $\Delta$ as a fan in $\R^{\mathcal{A}}/(1,\ldots,1)$.

   \begin{lemma}\label{lem-maxdim}
     Given a marked subdivision $T=\{(Q_l,\mathcal{A}_l)\}$ of
     $\Delta$ of type $\alpha$, we have
     \begin{displaymath}
       \dim(\alpha) \leq  \dim(C_T),
     \end{displaymath}
     where $C_T$ denotes the cone of the secondary fan corresponding to $T$. Equality holds if and only if in $T$ all lattice points of
     $\Delta$ are marked, i.e.\ if $\bigcup_{l} \mathcal{A}_l = \Delta\cap \Z^3$.
   \end{lemma}
The proof is analogous to Lemma 2.5 of \cite{MMS09}.

Since many tropical polynomials can induce the same tropical
 surface, the secondary fan is not the parameter space for
tropical surfaces. However, singular tropical surfaces of maximal-dimensional geometric type are parametrized by the union of (the interior of) codimension one cones of the secondary fan which correspond to dual marked subdivisions with all lattice points marked.
This feature also explains our interest in singular tropical surfaces of maximal-dimensional geometric type.

  \subsection{The tropicalisation of $\Sing_\bp(\Delta)=\ker(A)$}\label{subsec-family}

 We use the following known results about
   the tropicalisation of linear spaces (\cite{Stu02}, \S\;9.3, \cite{FS05}, \cite{AK06}). The
   tropicalisation of the linear space $\ker(A)$ depends only on the matroid $M$ associated to
   $A$ as follows: we define $M$ by its collection of circuits, which are minimal sets $\{i_1,\ldots,i_r\}\subset
   \{1,\ldots,s\}$ such that the columns $b_{i_1},\ldots,b_{i_r}$ of a
   Gale dual $B$ of $A$ are linearly dependent.
   A Gale dual is a matrix $B$ whose rows span the kernel of
   $A$.
   Given $u\in \R^s$, let $\mathcal{F}(u)$ denote the unique \emph{flag of subsets}
   \begin{displaymath}
     \emptyset=: F_0 \subsetneqq F_1 \subsetneqq \ldots \subsetneqq
     F_k\subsetneqq F_{k+1}:= \{1,\ldots,s\}
   \end{displaymath}
   such that
   \begin{displaymath}
     u_i<u_j\;\;\;\Longleftrightarrow\;\;\;\exists\;l\;:\; i\in F_{l-1}
     \mbox{ and } j\not\in F_{l-1}.
   \end{displaymath}
   In particular,
   \begin{displaymath}
     u_i=u_j\;\;\;\Longleftrightarrow\;\;\;\exists\;l\;:\; i,j\in
     F_l\setminus F_{l-1}.
   \end{displaymath}
   The \emph{weight class} of a flag $\mathcal{F}$ is the set of all
   $u$ such that $\mathcal{F}(u)=\mathcal{F}$.
  \longer{ We say that a flag of subsets $\mathcal{F}(u)$ is in the boundary
   of a flag of subsets $\mathcal{F}(v)$ if $\mathcal{F}(u)$ can be
   obtained from $\mathcal{F}(v)$ by removing some of the $F_i$.
   Note that the weight class of a flag is an open cone, and
   $\mathcal{F}(u)$ is in the boundary of $\mathcal{F}(v)$ if and only
   if $u$ is in the boundary of the weight class of $\mathcal{F}(v)$.}

   A flag $\mathcal{F}$ is a \emph{flag of flats} of the Gale dual $B$
   of $A$ respectively of the associated matroid $M$ if the linear span
   of the vectors $\{b_j \;|\; j\in F_i\}$ contains no $b_k$ with
   $k\notin F_i$. As before, the vectors $b_j$ denote the columns of $B$.
   It follows from Theorem 1 of \cite{AK06} resp.\
   Theorem 4.1 of \cite{FS05} that the Bergman fan of a matroid
   $M$ is the set of all weight classes of flags of flats of
   $M$.

   As a consequence, we can study $\Trop(\ker(A))$ by
   studying weight classes of flags of flats of a Gale dual of $A$.
   Note that since $A$ is a $4\times s$-matrix, maximal flags of flats can
   be identified with
     flags of $s-4$ subspaces $V_i\subset \R^{s-4}$:
     \begin{displaymath}
       \{0\} \subsetneqq V_1 \subsetneqq \ldots \subsetneqq V_{s-4},
     \end{displaymath}
     where each $V_i$ is generated by a subset of the column vectors
     $b_j$ of the Gale dual $B$ of $A$ indexed by the set $F_i$, and the vectors
     $\{b_j\;|\; j\in F_i\}$ are all the column vectors of the Gale dual
     that are contained in the subspace $V_i$.
     In particular, $F_{s-4}=\{1,\ldots,s\}$.
     We set $F_i':=F_i\setminus F_{i-1}$.
     Each $F_i'$ must of course consist of at least one element $j$.
     Since we have $s$ vectors in total, we have 4 ``extra'' vectors
     that can a priori belong to any of the $F_i'$.
     In the next lemma, we show how the four extra vectors can be spread.

   \begin{lemma}\label{lem-chains}
     With the notation from above, for each flag of
     flats $\mathcal{F}=\mathcal{F}(u)$ of a Gale dual $B$ of $A$ we have either
     \begin{enumerate}
     \item $\# F_i'=1$ for all $i=1,\ldots,s-5$ and $\# F_{s-4}'=5$, or
     \item $\# F_{s-4}'=4$ and there is a $j\in \{1,\ldots,s-5\}$ with $\# F_{j}'=2$, or
     \item $\# F_{s-4}'=3$ and there is a $j\in \{1,\ldots,s-5\}$ with $\# F_{j}'=3$, or
     \item $\# F_{s-4}'=3$ and there are $i<j\in \{1,\ldots,s-5\}$ with $\# F_{i}'=\#F_j'=2$.
     \end{enumerate}
     In each case, the lattice points corresponding to the indices in $ F_{s-4}'$ form a circuit.
     In the first case, this is a circuit of type (A) or (B) as in
     Remark \ref{rem-circuit}, in the second case of type (C) or (D),
     and in the third and fourth case of type (E).

     In the second case, all points $m_r$ with $r\in
     F_l'$, $l>j$, are on the same plane as the four points of $
     F_{s-4}'$, and none of the points with $r\in F_{j}'$ is on this
     plane.

     In the third case, all points $m_r$ with $r\in
     F_l'$, $l>j$, are on the same line as the three points of
     $F_{s-4}'$, and each choice of two of the points in $F_{j}'$
     spans the space together with the three points of
     $F_{s-4}'$.

     In the fourth case, all points $m_r$ with $r\in
     F_l'$, $l>j$, are on the same line as the three points of
     $F_{s-4}'$, and all points $m_r$ with $r\in
     F_l'$, $j>l>i$, are on the same plane as the three points of
     $F_{s-4}'$ and the two points of $F_j'$, and the two points of
     $F_i'$ do not lie on this plane.
   \end{lemma}
   The proof is a straightforward generalisation of Lemma 3.7 of
   \cite{MMS09}. Note that with this Lemma we describe only interior points of cones
   corresponding to weight classes of top dimension in $\Trop(\ker(A))$.
The analogous statement to Remark 3.8 of \cite{MMS09}
   holds true as well: for any circuit and any choice of points
   satisfying the affine dependencies as above we can find a
   corresponding weight class in $\Trop(\ker(A))$.  That means that whenever the coefficients of a tropical polynomial
   meet one of the above conditions, it lifts to a polynomial over $\K$ defining a surface with singularity at (1,1,1).

\section{The tropical discriminant revisited}\label{sec:tropicaldiscriminant}

For ${\bx}\in \R^n$ arbitrary, denote by $\bp_{\bx}\in(\K^\ast)^n$ a
point with $\val(\bp_{\bx})={\bx}$, and consider the family
$\Sing_{\bp_{\bx}}(\Delta)$ of surfaces with a singularity in
$\bp_{\bx}$. Its tropicalisation $\Trop(\Sing_{\bp_{\bx}}(\Delta))$
does not depend on the choice of $\bp_{\bx}$. Moreover, it follows
from Remark 3.2 of \cite{MMS09} that it is a shift of
$\Trop(\Sing_\bp(\Delta))=\Trop(\ker(A))$ by a vector which we denote
by $v({\bx})$ whose coordinates in $\R^s/(1,\ldots,1)$ are given by the scalar
products of the $m\in \mathcal{A}$ with ${\bx}$.

If we let ${\bx}$ vary over all points in $\R^n$, it follows that
$v({\bx})$ varies over all points in the rowspace of the matrix $A$ in
$\R^s/(1,\ldots,1)$. In the following, we denote the rowspace of $A$
in $\R^s/(1,\ldots,1)$ by $L$. Notice that $L$ also equals the lineality space of the
secondary fan.

\begin{notation}\label{not-v}
  Let $v:\R^n\rightarrow L$ denote the linear map sending ${\bx}$ to
  $v({\bx})=(m\cdot {\bx})_{m\in\mathcal{A}}$ as above. Notice that $v$
  is a bijective linear map between vector spaces of dimension $n$.
\end{notation}

This illustrates the equality $\Trop(\ker(A))+\rowspace(A)=\Trop(\Sing(\Delta))$ which is proved in Theorem 1.1 of \cite{DFS05}.
Since we assume that $\Delta$ yields a non-defective point configuration, it follows from \cite{GKZ}, 10.1.2, that
$\Trop(\Sing(\Delta))$ is a subfan of the codimension-one-skeleton of the secondary fan. Therefore it comes with a
natural fan structure given by the secondary fan.
Since it equals $\Trop(\ker(A))+\rowspace(A)$, it also comes with a natural fan structure by weight classes of the lattice
of flats of the matroid of $A$.
In general, these two fan structures are not compatible --- cones can overlap, cut through cones, be smashed to lower
dimension etc.\ In the following, we define the notion of a generic tropical surface and restrict our results in
Theorems \ref{thm1} and \ref{thm2} to generic surfaces --- these are surfaces for which the two fan structures locally
around the coefficient vector $u$ are best compatible. The set of generic surfaces is of top dimension.

\begin{notation}\label{not-C(a)}
 For a point $u\in \R^s/(1,\ldots,1)$, we set $C(u)$ the unique cone
 of the secondary fan with $u\in \mbox{relint}(C(u))$.
 Notice that $C(u)=C(u+l)$ for every $l\in L$.
\end{notation}

\begin{definition}\label{def:defective}
  We call a weight class $C$, i.e.\ a cone of $\Trop(\ker(A))$,
  \emph{defective} if there exists a point $u\in C+L$ with
  $\dim(C+L)<\dim(C(u))$.
\end{definition}

\begin{remark}
  If $C$ is a weight class and $u\in C$ such that $C(u)$ has
  codimension one in the secondary fan, then $C$ is defective if and
  only if $\Span(C)\cap L\not=\{0\}$.
\end{remark}

\begin{example}
  We consider the point configuration
  $\mathcal{A}=\{m_a,m_b,\ldots,m_h\}$ with
  \begin{displaymath}
    \begin{array}{llll}
      m_a=(0,0,0), & m_b=(0,0,1), & m_c=(0,0,2), & m_d=(0,1,0),\\
      m_e=(0,-1,0),& m_f=(1,0,0), & m_g=(1,1,0), & m_h=(-1,0,0)
    \end{array}
  \end{displaymath}
  and we consider the weight class
  \begin{displaymath}
    C=\{x_{m_a}=x_{m_b}=x_{m_c}>x_{m_d}=x_{m_e}>x_{m_f}=x_{m_g}>x_{m_h}\}.
  \end{displaymath}
  The corresponding subdivision of the polytope $\Delta$ is shown in
  Figure \ref{fig:defectiveweightclass}.
  \begin{figure}[h]
    \centering
    \includegraphics[width=8cm]{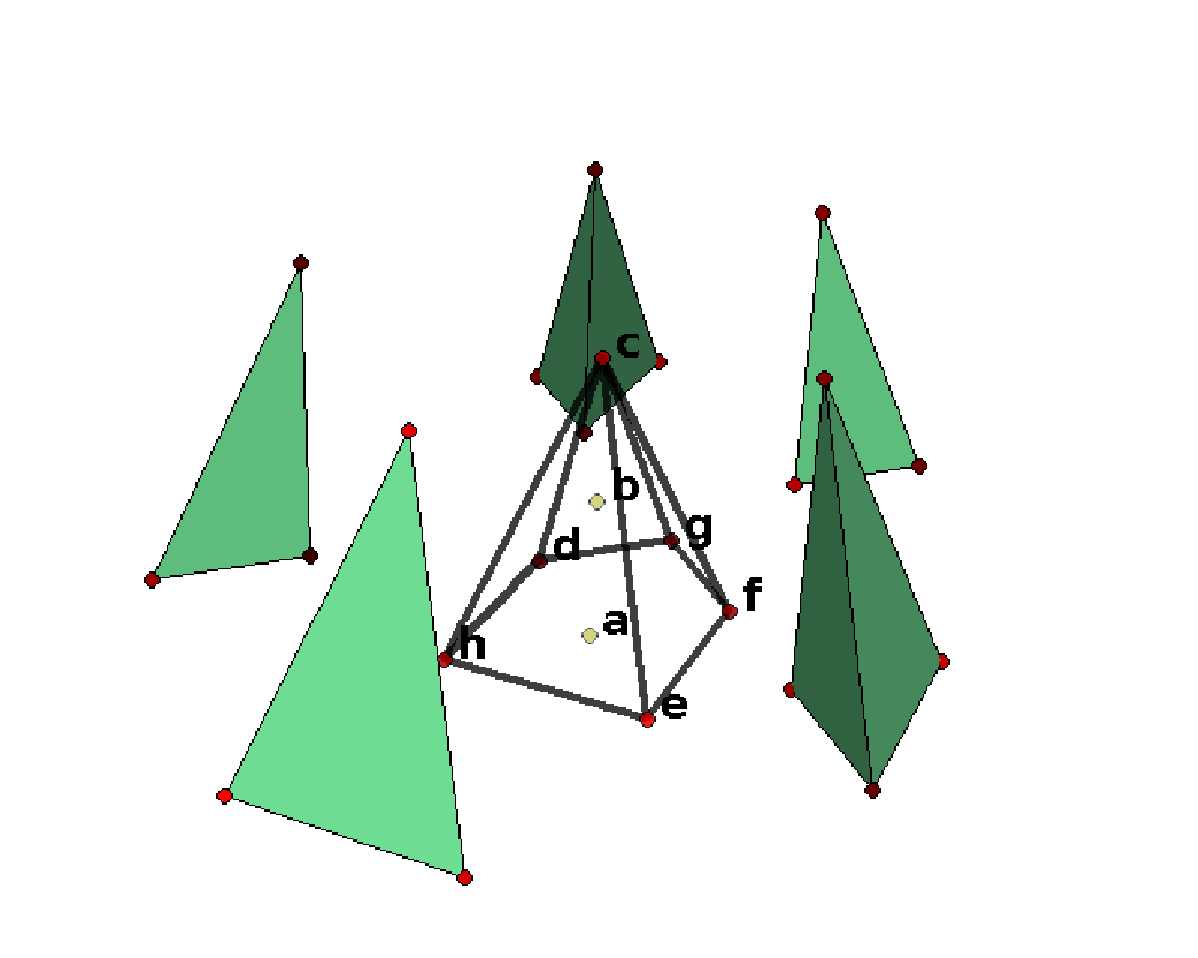}
    \caption{A subdivision corresponding to a defective weight class.}
    \label{fig:defectiveweightclass}
  \end{figure}
 For a point $u$ in the weight class $C$, the
  corresponding cone $C(u)$ in the secondary fan is of codimension
  one. However, the intersection of $\Span(C)$ with the lineality
  space in $\R^8/(1,\ldots,1)$ is $1$-dimensional, since it contains
  the vector $(0,0,0,0,0,-1,-1,1)$. This shows that the weight class
  is defective.

  Indeed, the weigth class $C$ shares a facet with each of the two
  weight classes
  \begin{displaymath}
    C'=\{x_{m_a}=x_{m_b}=x_{m_c}>x_{m_d}=x_{m_e}>x_{m_f}=x_{m_h}>x_{m_g}\}
  \end{displaymath}
  and
  \begin{displaymath}
    C''=\{x_{m_a}=x_{m_b}=x_{m_c}>x_{m_d}=x_{m_e}>x_{m_g}=x_{m_h}>x_{m_f}\}.
  \end{displaymath}
  The span of each of these two weight classes intersects the
  lineality space transversally. The cone $C(u)$ from above is
  just the union
  \begin{displaymath}
    (C+L)\cup(C'+L)\cup (C''+L),
  \end{displaymath}
  where actually $C+L$ is not needed, since it is a face of both
  $C'+L$ and $C''+L$. This is thus an example that a full-dimensional
  weight class in $\Trop(\Ker(A))$ may lead to a lower dimensional
  cone in the tropical discriminant of $A$ which lies in the interior
  of a full-dimensional cone of the tropical discriminant.

  Note that in this example the point configuration $\mathcal{A}$ itself is
  not defective, however the subset consisting of points $m_a,m_b,m_c,m_d,m_e,m_f $ is.
\end{example}

Assume $C$ is a non-defective weight class, then $C+L$ is contained in
cones of the secondary fan of dimension equal to $\dim(C+L)$ or
less. The set of all $u\in C$ with $\dim(C+L)>\dim(C(u))$ is obviously
of smaller dimension than $\dim C$.

We now define the notion of a generic tropical surface. We choose
this definition in such a way that if we
decompose the coefficient
vector $u$ of a generic surface as $v+l$ with $v\in \Trop(\ker(A))$
and $l\in L$, and $C$ is the weight class of $\Trop(\ker(A))$
containing $v$ in its relative interior, then all dimensions are as
expected, i.e.\ $C$ is of top dimension and $\dim(C+L)=\dim C(u)$.

\begin{definition}\label{def-generic1}
  We call a point $u\in \Trop(\Ker(A))+L\subseteq\R^s/(1,\ldots,1)$ in the tropical discriminant of
  $\mathcal{A}$ \emph{generic} if it lies outside the locus formed
  by $C+L$, where a cone $C$ of $\Trop(\Ker(A))$ either is
  defective, or is not of the top dimension, or satisfies
  $\dim(C+L)>\dim(C(v))$.
 The singular
  tropical hypersurface defined by the tropical polynomial $F_u$ is
  then also called \emph{generic}.
\end{definition}

From the above, it is obvious that the set of generic points in the tropical discriminant is of top dimension.

Note that in Theorem \ref{thm2} we consider generic surfaces whose
dual marked subdivision is of codimension one in the secondary fan.
For defective point configurations, such surfaces do not exist.

\begin{proof}[Proof of Theorem \ref{thm1}:]
Let $u\in \Trop(\Sing(\Delta))$ be generic. It follows from the
definition of genericity that we can write
$u$ as a sum $v+l$ with
$v\in \Trop(\ker(A))$ and $l\in L$, such that the weight class $C$
of $\Trop(\ker(A))$ which contains $v$ in its relative interior is
top-dimensional and satisfies $\dim(C+L)=\dim(C(v))$. Assume
$C(v)=C(u)$ is a cone of codimension $c$ of the secondary fan.
Notice that the representation of $u$ as a sum as above is not
unique. Firstly, there might be several weight classes $C$ in
$\Trop(\ker(A))$ such that we can write $u$ as the sum of a vector
in $C$ and a vector in $L$. Secondly, even if we fix one cone $C$,
there might be several representations of $u$ as the sum of a vector
in this $C$ and a vector in $L$. For now, let us fix one weight
class $C$ which allows a representation of $u$ as $u=v+l$ with $v\in
C$ and $l\in L$.


Since $\dim \Trop(\ker(A))=s-1-(n+1)$ (where $s=\#\mathcal{A}$) and
$v\in C$ is in a
top-dimensional weight class, we have
$\dim(C+L)=\dim(C)+\dim(L)-\dim(\Span(C)\cap
L)=s-1-(n+1)+n-\dim(\Span(C)\cap L) = \dim(C(v))=s-1-c$, where
$\Span(C)$ denotes the smallest linear space containing $C$. It
follows that $\dim(\Span(C)\cap L)= c-1$. Therefore there exists a
$c-1$-dimensional polyhedron in $H\subset C$ such that for all $h\in
H$ we have $v+h\in C$. We can thus write $v$ also as $v=(v+h)-h$,
where the first summand is in $C$ and the second summand is in $L$,
and these are all possibilities to represent $v$ as a sum of a
vector in $C$ plus a vector in $L$. Consequently, we can write $u$
as $u=(v+h)+(l-h)$ and again, these are all possibilities to
represent $u$ as a sum of a vector in $C$ and a vector in $L$. It
follows that $F_{u}$ defines a tropical surface which is singular at
all points $\bx_{l-h}$, where $\bx_{l-h}\in \R^n$ denotes the
preimage of the bijective linear map sending ${\bx}\in \R^n$ to $
v({\bx})=(m\cdot {\bx})_{m\in \mathcal{A}}$ from Notation
\ref{not-v}. Since the map $v^{-1}$ maps the $c-1$-dimensional
polyhedron $l-H$ to a $c-1$-dimensional polyhedron, it follows that
all singular points of the surface of $F_{u}$ that we get by
decomposing $u$ as a sum of a vector in $C$ and a vector in $L$ lie
in a $c-1$-dimensional polyhedron.
 As we have seen above there may be several (but finitely many) weight
 classes $C$ in $\Trop(\ker(A))$ such that we can write $u$ as the sum
 of a vector in $C$ and a vector in $L$, and it thus follows that the set
 of singular points of the tropical surface defined by $F_{u}$ is a finite union of $c-1$-dimensional polyhedra.
\end{proof}

\begin{remark}\label{rem:ex-thomas}
  Recall again Example \ref{ex-thomas} where we had a surface $S$ with two
  singular points. These two singular points arise because we can interpret the coefficient
  vector $u$ of the tropical polynomial defining $S$ in two ways as a sum of a vector in a
  weight class of $\Trop(\ker(A))$ and a vector in the lineality space. The two weight classes are different.
 The point configuration in question corresponds to
  the matrix
  \begin{displaymath}
    A=\left(\begin{array}{cccrccc}1&1&1&1&1&1&1\\0&0&0&-1&0&1&2\\0&0&0&-1&1&0&1\\0&1&2&0&0&0&1\end{array}\right),
  \end{displaymath}
  and the singular point $G=(0,0,0)$ on $F_u$ with
  \begin{displaymath}
    u=(0,0,0,-8,-5,-5,-5)\in\Trop(\Ker(A))
  \end{displaymath}
  comes from the weight class containing $u$. However, we can also
  write $u$ as
  \begin{displaymath}
    u=(0,0,0,-6,-6,-6,-8)+(0,0,0,-2,1,1,3)=v+l
  \end{displaymath}
  where
  \begin{displaymath}
    l=(0,0,0,-1,0,1,2)+(0,0,0,-1,1,0,1)=(0,0,0,-2,1,1,3)
  \end{displaymath}
  belongs to the lineality space of the secondary fan of the point
  configuration and $v$ belongs to some other
  weight class. The corresponding singular point on $S$ is
  $H=(-1,-1,0)$, since we have added once the vector of
  $x$-coordinates and once the vector of $y$-coordinates to the weight
  vector $v$ in the weight class in order to get $u$. This corresponds
  to shifting the whole surface (determined by $F_v$, which is singular at $0$) by
  $(-1,-1,0)$ (see also \ref{not-v} and before).
  Examples \ref{ex:barycenter1} and \ref{ex:barycenter2} give
  further explanations concerning this example.

  This shows that even if the point $u$ in the tropical discriminant
  is generic, the surface corresponding to $u$ may have more than one
  singular point.
\end{remark}

\section{The classification}\label{sec-class}
Now, using the preparation from Section \ref{sec-basic}, we prove Theorem \ref{thm2}. In particular, we consider the points $u\in \Trop(\Sing(\Delta))$
which are generic in the sense of Definition \ref{def-generic1}, and
which in addition satisfy $\dim(C(u))=s-2$, where $C(u)$ is as in
Notation \ref{not-C(a)}. In addition, we work in the situation where
the dual marked subdivision as in Subsection
\ref{subsec-hypersurf} has all lattice points marked (see Lemma
\ref{lem-maxdim}).
Since we can always  write $u=v+l$ for some $a\in \Trop(\ker(A))$ and
$l\in \rowspace(A)$, just as in the proof of Theorem \ref{thm1}
above, we can classify the singularities of the tropical surface
defined by $F_v$ with $v\in \Trop(\ker(A))$ first, and then
investigate how the shift to $F_{u}$ effects the location of the
singular points.
We thus have to consider all different types of weight classes as in
Lemma \ref{lem-chains}, and the corresponding possible types of
circuits.
It turns out that in most cases we do not have to worry about the
shift when passing from $F_v$ to $F_{u}$, since we describe the
location of the singular point relative to other points in the
surface, e.g.\ as the midpoint of an edge. This midpoint is of course
shifted accordingly.

\subsection{Weight class as in Lemma \ref{lem-chains}(a), circuit (A) of Remark \ref{rem-circuit}}
Let $u\in \Trop(\ker(A))$ be in a weight class as in Lemma
\ref{lem-chains}(a), and assume $F_{s-4}'=\{a,b,c,d,e\}$.
Consider
the marked subdivision defined by $u$ as in Subsection
\ref{subsec-hypersurf}. As the heights of the points
$m_a,m_b,m_c,m_d$ and $m_e$ are biggest, it follows that the convex
hull spanned by these points is a polytope of the subdivision. Let
us first assume that this polytope is a circuit of type (A) as in
Remark \ref{rem-circuit}. The vertex of the tropical surface dual to
this pentatope is at the point
 $(x,y,z)$ where the maximum is
     attained by the corresponding five terms of $\trop\{u_m+m\cdot (x,y,z)\}$, in particular the five terms are
     equal at this vertex. That means, we can set the five terms equal
     and solve for $x,y$ and $z$ to get the position of the
     vertex. But since the coefficients are all equal, we get $x=y=z=0$
     when solving. Notice that $(0,0,0)$ is the tropicalisation of the singular point $(1,1,1)$.

Since we require that all lattice points are marked, this polytope
cannot contain any lattice point besides these five. By Theorem 3.5
of \cite{Sca85}, a pentatope which does not contain any lattice
point besides its five vertices are IUA-equivalent to the tuple of
points $(0,0,0)$, $(1,0,0)$, $(0,1,0)$, $(0,0,1)$ and $(1,p,q)$ with $p$ and $q$ coprime. It is a bipyramid.

It follows that in this situation the node of the tropical surface is at a vertex with six adjacent edges and nine adjacent $2$-dimensional polyhedra.

This settles case (a.1) of Theorem \ref{thm2}.
\subsection{Weight class as in Lemma \ref{lem-chains}(a), circuit (B) of Remark \ref{rem-circuit}}

As above, it follows that the singular point $(0,0,0)$ is dual to
the convex hull of $m_a,m_b,m_c,m_d$ and $m_e$. This
is a vertex of
the tropical surface with four adjacent edges and six
$2$-dimensional polyhedra, just as a smooth vertex. However, if we
define the multiplicity of a vertex of a tropical hypersurface
analogously to the case of tropical curves as the lattice volume of
the corresponding polytope in the dual subdivision, then it follows
that the singular point is a vertex of higher multiplicity. More
precisely, the multiplicity can be 4,5,7,11,13,17,19 or 20. This
follows from the classification of $3$-dimensional tetrahedra with
one interior lattice point (and no other lattice points besides the
vertices) (see \cite{Rez06}, Theorem 7). Since we require that all
lattice points are marked, the tetrahedron which is the convex hull
of $m_a,m_b,m_c,m_d$ and $m_e$ has to be of this form. The
classification states that such a tetrahedron is IUA-equivalent to
one of the following eight: it has vertices $(0,0,0)$, $(1,0,0)$,
$(0,1,0)$ and, respectively, $(3,3,4)$, $(2,2,5)$, $(2,4,7)$,
$(2,6,11)$, $(2,7,13)$, $(2,9,17)$, $(2,13,19)$, or $(3,7,20)$. This
settles case (a.2) of Theorem \ref{thm2}.

\subsection{Weight class as in Lemma \ref{lem-chains}(b), circuit (C) of Remark \ref{rem-circuit}}
Let $F_{s-4}'=\{a,b,c,d\}$ and $F_j'=\{e,f\}$.
\subsubsection {Assume that in the subdivision, $m_e$ and $m_f$ both form a pyramide with the triangle spanned by $m_a,m_b,m_c$ and $m_d$ as base.}\label{subsubsec:b.1.1}
In particular, $m_e$ and $m_f$ must lie on different sides of the plane spanned by $m_a$, $m_b$, $m_c$ and $m_d$. 
Since there are no other circuits, and we require that all lattice points in
$\Delta$ are marked, both of these pyramids contain no further lattice points.
\begin{lemma}\label{lem-13}
 Let four lattice points $m_a,m_b,m_c$ and $m_d$ in an affine plane in $\R^3$ form a circuit of type (C) as in Remark
 \ref{rem-circuit}. Let $m_e$ be a fifth lattice point that forms a pyramid with this circuit as base and assume this pyramid
 contains no further lattice points. Then $m_e$ has integral distance $1$ or $3$ from the plane spanned by the circuit.
\end{lemma}
\begin{proof}
 We can assume that the plane spanned by $m_a,m_b,m_c$ and $m_d$ is the $x=0$-plane, and,
 using a suitable automorphism of $\Z^2$, we can bring
 these four points to $(0,0)$, $(-1,-1)$, $(-2,-1)$ and $(-1,-2)$. Denote the triangle spanned by these points by $T$. Also we
 assume without restriction that the $x$-coordinate of $m_e$ is positive. We have to show that it is then either $1$ or $3$.
Consider a lattice point $m$ with $x$-coordinate $1$ and let $C_m$ be the cone with vertex $m$ and spanned by the rays $m$,
$m-(0,-2,-1)$ and $m-(0,-1,-2)$. Intersect this cone with the plane $x=k$ for some choice of $k>1$ (see Figure \ref{fig:conem}).
 \begin{figure}[h]
    \centering
\begin{picture}(0,0)%
\includegraphics{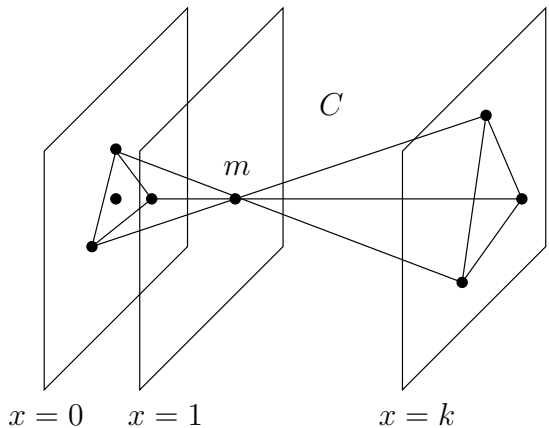}%
\end{picture}%
\setlength{\unitlength}{3947sp}%
\begingroup\makeatletter\ifx\SetFigFont\undefined%
\gdef\SetFigFont#1#2#3#4#5{%
  \reset@font\fontsize{#1}{#2pt}%
  \fontfamily{#3}\fontseries{#4}\fontshape{#5}%
  \selectfont}%
\fi\endgroup%
\begin{picture}(3402,2706)(3811,-4255)
\put(5776,-2236){\makebox(0,0)[lb]{\smash{{\SetFigFont{12}{14.4}{\familydefault}{\mddefault}{\updefault}{\color[rgb]{0,0,0}$C$}%
}}}}
\put(3826,-4186){\makebox(0,0)[lb]{\smash{{\SetFigFont{12}{14.4}{\familydefault}{\mddefault}{\updefault}{\color[rgb]{0,0,0}$x=0$}%
}}}}
\put(4576,-4186){\makebox(0,0)[lb]{\smash{{\SetFigFont{12}{14.4}{\familydefault}{\mddefault}{\updefault}{\color[rgb]{0,0,0}$x=1$}%
}}}}
\put(6151,-4186){\makebox(0,0)[lb]{\smash{{\SetFigFont{12}{14.4}{\familydefault}{\mddefault}{\updefault}{\color[rgb]{0,0,0}$x=k$}%
}}}}
\put(5176,-2611){\makebox(0,0)[lb]{\smash{{\SetFigFont{12}{14.4}{\familydefault}{\mddefault}{\updefault}{\color[rgb]{0,0,0}$m$}%
}}}}
\end{picture}%

    \caption{The cone $C_m$.}
    \label{fig:conem}
  \end{figure}

For any lattice point in $C_m\cap \{x=k\}\cap \Z^3$, consider the
pyramid that this point forms with $T$ as base. This pyramid will
contain the point $m$. If we move $m$ by a step of integer length
$1$, the triangle $C_m\cap \{x=k\}$ is shifted by $k$. Compared to
$T$ the triangle $C_m\cap \{x=k\}$ is grown by a factor of $k-1$. Of
course, we can also move $m$ to a point with a different
$x$-coordinate, this will add more triangles (smaller in size) such
that for each point inside a triangle we know that the corresponding
pyramid contains another lattice point. We show that for $k\neq 3$
the
shifted triangles cover all lattice points with $x$-coordinate
$k$. It follows that any pyramid with $T$ as base and with a vertex
with $x$-coordinate $k\neq 1,3$ contains another lattice point.
Figure \ref{fig:shifts} shows the plane $\{x=k\}$ with the
$k$-shifts of the triangle $C_m\cap \{x=k\}$.
\begin{figure}[h]
 \centering
\begin{picture}(0,0)%
\includegraphics{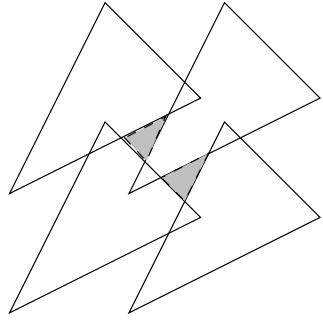}%
\end{picture}%
\setlength{\unitlength}{3947sp}%
\begingroup\makeatletter\ifx\SetFigFont\undefined%
\gdef\SetFigFont#1#2#3#4#5{%
  \reset@font\fontsize{#1}{#2pt}%
  \fontfamily{#3}\fontseries{#4}\fontshape{#5}%
  \selectfont}%
\fi\endgroup%
\begin{picture}(1974,1974)(4189,-4273)
\end{picture}%

\caption{The shifts of the triangle $C_m\cap \{x=k\}$ on the plane $\{x=k\}$.}
\label{fig:shifts}
\end{figure}

 Let us compute the vertices of the right shaded region which is not yet covered by a triangle.
 Assume the left most vertex of the top right triangle has coordinates $(1,1)$ in the plane, then
 it follows that the coordinates of the vertices of the shaded region are $(\frac{2k}{3}-1,\frac{k}{3})$,
 $(\frac{2k}{3}+1,\frac{k}{3}+1)$ and $(\frac{2k}{3},\frac{k}{3}-1)$. Independently of $k$, this is a
 triangle of lattice area $3$, with the point $(\frac{2k}{3},\frac{k}{3})$ as midpoint from which we
 reach the three vertices by a lattice step to the left, down, and to the upper right.
This triangle has an interior lattice point if and only if $k$ is divisible by $3$. In this case, the
lattice point is $(\frac{2k}{3},\frac{k}{3})$ (see Figure \ref{fig:kregion}).
\begin{figure}[h]
 \centering
\begin{picture}(0,0)%
\includegraphics{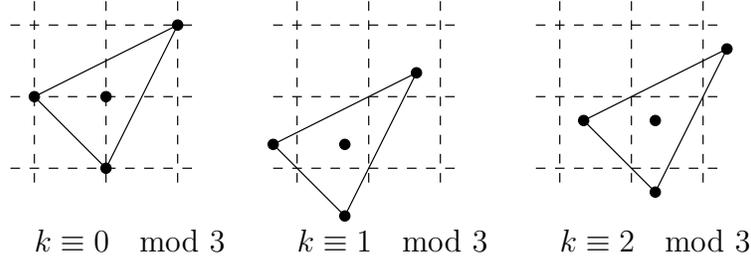}%
\end{picture}%
\setlength{\unitlength}{3947sp}%
\begingroup\makeatletter\ifx\SetFigFont\undefined%
\gdef\SetFigFont#1#2#3#4#5{%
  \reset@font\fontsize{#1}{#2pt}%
  \fontfamily{#3}\fontseries{#4}\fontshape{#5}%
  \selectfont}%
\fi\endgroup%
\begin{picture}(4552,1665)(9139,-3814)
\put(9301,-3736){\makebox(0,0)[lb]{\smash{{\SetFigFont{12}{14.4}{\familydefault}{\mddefault}{\updefault}{\color[rgb]{0,0,0}$k\equiv 0\; \mod 3$}%
}}}}
\put(10951,-3736){\makebox(0,0)[lb]{\smash{{\SetFigFont{12}{14.4}{\familydefault}{\mddefault}{\updefault}{\color[rgb]{0,0,0}$k\equiv 1\; \mod 3$}%
}}}}
\put(12601,-3736){\makebox(0,0)[lb]{\smash{{\SetFigFont{12}{14.4}{\familydefault}{\mddefault}{\updefault}{\color[rgb]{0,0,0}$k\equiv 2\; \mod 3$}%
}}}}
\end{picture}%

\caption{The non-covered region for different values of $k$.}
\label{fig:kregion}
\end{figure}

Analogously, we can compute the vertices of the left shaded region
and see that it has an interior lattice point if and only if $k$ is
divisible by $3$, and then this lattice point has coordinates
$(\frac{k}{3},\frac{2k}{3})$. It follows that for any $k$ which is
not divisible by $3$ the $k$-shifts of the triangle $C_m\cap
\{x=k\}\cap \Z^3$ cover already all lattice points. That is, any
pyramid with $T$ as base and with a vertex with $x$-coordinate
which
is not divisible by $3$ contains another lattice point with
$x$-coordinate $1$.

If $k$ is divisible by $3$, then a pyramid with a vertex with
coordinates $(k,\frac{2k}{3}+ik,\frac{k}{3}+jk)$
or
$(k,\frac{k}{3}+ik,\frac{2k}{3}+jk)$ where $i,j\in \Z$ does not
contain a lattice point with $x$-coordinate $1$. Here, we take the
effect of the $k$-shifts of $C_m\cap \{x=k\}\cap \Z^3$ on the shaded
regions into account. Let us now assume that $k=3^l\cdot h$, where
$h\neq 1$ and $3\nmid h$. Now move $m$ to a point with
$x$-coordinate $3^l$. It follows using the same arguments as above
that any pyramid with base $T$ and a vertex with $x$-coordinate $k$
contains a lattice point with $x$-coordinate $3^l$. Next assume that
$k=3^l$, $l\geq 2$. Using $m$ with $x$-coordinate $1$ as before we
see that the only possibilities to get a pyramid which does contain
a lattice point with $x$-coordinate $1$ are that the vertex has
$(y,z)$-coordinates divisible by $3^{l-1}$ and not by $3^l$. Using
$m$ with $x$-coordinate $3$ we see that the only possibilities to
get a pyramid which does not contain a lattice point with
$x$-coordinate $3$ are that the vertex has $(y,z)$-coordinates
divisible by $3^{l-2}$ and not by $3^{l-1}$. As there is no vertex
which satisfies both it follows that any pyramid with a vertex with
$x$-coordinate $k=3^l$, $l\geq 2$, contains a lattice point with
$x$-coordinate $1$, or it contains a lattice point  with
$x$-coordinate $3$. In any case, it contains another lattice point.
It follows that the $x$-coordinate of the vertex $m_e$ can only be
$1$ or $3$.
\end{proof}
\begin{remark}
 There are vertices $m_e$ with integral distance $1$ and $3$ to the plane containing the circuit of type (C) of
 Remark \ref{rem-circuit} such that the pyramid formed by the circuit and $m_e$ contains no further lattice points,
 e.g. the convex hull of the points $(0,1,0)$, $(0,0,1)$, $(0,2,2)$ and $(3,0,2)$, or the convex hull of the points
 $(0,0,0)$, $(0,1,2)$, $(0,2,1)$ and $(1,0,0)$.
\end{remark}

Now solve the equations given by the tropical polynomial to get the positions of the two vertices corresponding to the two pyramids.
The $x$-coordinates of $m_e$ and $m_f$ can either be the negative of each other, or one can be $3$ and the other $-1$.
Since $m_a,m_b,m_c$ and $m_d$ have biggest and equal height, it
follows that the edge dual to the convex hull of $m_a,m_b,m_c$ and
$m_d$ satisfies the equations $y=0$ and $z=0$. If $\lambda=u_{m_a}$ is
the biggest weight (the weight of $m_a,m_b,m_c$ and $m_d$), and
$\mu=u_{m_e}$ is the weight of $m_e$ and $m_f$, it follows that the
vertex dual to the pyramid with vertex $m_e$ is at $(\mu-\lambda,0,0)$
(resp.\ $(\frac{1}{3}\cdot (\mu-\lambda),0,0)$) and the vertex dual to
the pyramid with vertex $m_f$ is at $(\lambda-\mu,0,0)$ (resp.\
$(\frac{1}{3}\cdot (\lambda-\mu),0,0)$). It follows that the singular
point $(0,0,0)$ is either exactly in the middle of the edge dual to
the convex hull of $m_a,m_b,m_c$ and $m_d$, or subdivides the edge
into parts whose distances have ratio 1:3. This explains the first cases of
(b.1) in Theorem \ref{thm2}.

\subsubsection {Assume that at most one of the points $m_e$ and $m_f$ forms a pyramide with the triangle spanned by $m_a$, $m_b$, $m_c$ and $m_d$ as base.} \label{subsubsec-eq}

As before, assume that $m_a=(0,0,0)$, $m_b=(0,1,1)$, $m_c=(0,2,1)$
and $m_d=(0,1,2)$. 
Any point $m$ that forms a pyramide with the triangle as base must have the absolute value of the $x$-coordinate $1$ or $3$ due to Lemma \ref{lem-13},
since otherwise we would
get extra lattice points, contradicting our
assumption that the plane is of maximal-dimensional geometric type. 
If $m$ is a point different from $m_e$ and $m_f$ but with $x$-coordinate $1$ or $-1$, then it cannot form a pyramide with the triangle as base since its coefficient is too low.
Thus we can conclude without restriction that in this situation, there is a pyramide with the triangle as base and with a vertex $m$ with $x$-coordinate $3$. 

Assume that $m_e$ does not form a pyramide with the triangle as base. We now determine the possible $x$-coordinates of $m_e$. We have seen already that then there is a pyramide with vertex $m$ with $x$-coordinate $ 3$. Because $\mu=u_{m_e}$ is the second biggest coefficient, the $x$-coordinate of $m_e$ must be  smaller than $3$. It cannot be $2$ however, since then by Lemma \ref{lem-13} the pyramide formed by $m_e$ and the triangle contains further lattice points. Even if this pyramide is not part of the subdivision, these additional lattice points would be contained in the convex body spanned by $m$, $m_e$ and the triangle. Because of their coefficients smaller $\mu$, they cannot be marked points of the subdivision, contradicting our assumption that the surface is of maximal-dimensional geometric type. With the same arguments, $m_e$ cannot have $x$-coordinate $-2$ or smaller $-3$. It follows that it must have the absolute value of the $x$-coordinate $1$.

There are a priori several possibilities for weight classes as in Lemma \ref{lem-chains}(b) from which our subdivision can arise.
In order to determine these possibilities, we have to decompose the coefficient vector $u\in
\Trop(\Sing(\Delta))$ of our tropical polynomial as a sum $v+l$ where $v$ is in a feasible weight class and $l$ is in the rowspace of $A$. Assume we have already added vectors of the rowspace to $u$ to achieve that the four points of the circuit have equal and biggest coefficients. Next we add a multiple of the vector of $x$-coordinates to make two coefficients of points outside the plane of the circuit equal and second biggest, the two points $m_e$ and $m_f$. To all points with $x$-coordinate one, we add the same value by adding the multiple of this rowspace vector. Thus there is a unique point with $x$-coordinate $\pm 1$ which is a candidate to be $m_e$ resp.\ $m_f$ --- the one with the biggest coefficient after adding rowspace vectors that make the coefficients of the circuit equal.
Candidates for $m_e$ and $m_f$ are now points with $x$-coordinate $\pm 3$ that form a pyramide with the triangle as base, and points with $x$-coordinate $\pm 1$ whose coefficient is biggest after adding rowspace vectors to make the coefficients of the circuit equal. Also, $m_e$ and $m_f$ must have different $x$-coordinates since otherwise the weight class would not intersect the corresponding weight class transversely which contradicts our assumtion that $u$ is generic (see Definition \ref{def-generic1}).

%
We therefore have the following four possibilities for weight classes (without restriction):
\begin{itemize}
 \item Let $m$ with $x$-coordinate $3$ form a pyramide with the triangle, and let $m_e$ be a point with $x$-coordinate one. Let $m_f$ with $x$-coordinate $-1$ form a pyramide with the triangle.
\item Let $m$ with $x$-coordinate $3$ form a pyramide with the triangle, and let $m_e$ be a point with $x$-coordinate one. Let $m_f$ with $x$-coordinate $-3$ form a pyramide with the triangle.
\item Let $m_1$ with $x$-coordinate $3$ and $m_2$ with $x$-coordinate $-3$ form a pyramide with the triangle. Let $m_e$ be a point with $x$-coordinate $1$ and $m_f$ with $x$-coordinate $-1$.
\item Let $m=m_f$ with $x$-coordinate $3$ form a pyramide with the triangle, and let $m_e$ be a point with $x$-coordinate one.
\end{itemize}
In each of the four cases, the rowspace of $A$ intersects the corresponding weight class transversely, and so there is at most one solution to decompose $u$ as a sum. The decomposition must be possible in at least one of the cases. 

The fourth case has to be treated separately. Note that the fourth case is the only one which can also arise if the edge dual to the triangle is unbounded.

In the first three cases, we introduce the notion of a \emph{virtual edge} dual to the triangle. This virtual edge is just as the actual edge dual to the triangle contained in the line $y=z=0$, however it ends at points whose $x$-coordinates differ from the actual $x$-coordinates of the vertices dual to the pyramides adjacent to the triangle. For a fixed weight class, i.e.\ for a fixed choice of $m_e$ and $m_f$ as above, we define the \emph{virtual vertex} corresponding to $m_e$ to be the vertex dual to the pyramide formed by the triangle and $m_e$ (even though this pyramide is not part of the subdivision). In the third case, we also define the virtual vertex of $m_f$ analogously. The virtual edge connects the virtual vertex of $m_e$ with the (virtual or actual) vertex of $m_f$. Note that the virtual edge contains the actual edge.
It follows from the previous Subsection that in the first and third case, the singular point is the midpoint of the virtual edge while in the second case, it subdivides the virtual edge with ration $1:3$.

Let us treat the first case exemplarily with more details. Denote by $\lambda$ the coefficient of $m_a$, $m_b$, $m_c$ and $m_d$, by $\mu$ the coefficient of $m_e$ and $m_f$ and by $\nu$ the coefficient of the point $m$ which forms a pyramide with the triangle. We have $\nu<\mu<\lambda$. 
The virtual vertex of $m_e$ has coordinates $(\lambda-\mu,0,0)$, the actual vertex --- i.e.\ the vertex corresponding to the pyramide formed by $m$ and the triangle --- has coordinates $(\frac{\lambda-\nu}{3},0,0)$.
Since $m_e$ has $x$-coordinate $1$ and $m$ has $x$-coordinate $3$ but forms a pyramide with the triangle, we must have $\mu<\frac{2\lambda+\nu}{3}$. This shows that the virtual edge is indeed longer than $E$.

In the fourth case, we cannot describe the location of the singular
point as some sort of midpoint as in the earlier cases,
a description
which does not change when we shift. When we solve for the position
of $V$ as before, and denote by $\lambda=u_{m_a}$ the highest
weight, i.e.\ the coefficient of the four points $m_a,\ldots,m_d$,
and by $\mu=u_{m_e}$ the coefficient of $m_e$ and $m_f$, then as
before we get $(\frac{1}{3}\cdot (\lambda-\mu),0,0)$ for the
coordinates of $V$. The singular point is at $(0,0,0)$ which is a
point of distance $\frac{\lambda-\mu}{3}$ from $V$. This distance
will not change of course when we shift, however the coefficients
$\lambda$ and $\mu$ are going to be changed by adding a vector in
the rowspace of $A$. Since there is a unique way of writing
$u$ as a sum of a vector in the weight class and a vector in the rowspace, we can in fact solve for the vector in the
rowspace which we need.
By our choice of coordinates for the point configuration, we can
deduce that we need to add the vector of $y$-coordinates in the
rowspace
$(u_{m_b}-u_{m_c})$-times and the vector of
$z$-coordinates
$(u_{m_b}-u_{m_d})$-times.
Then the four new coefficients of the circuit are equal, we have
\begin{align*}
  \lambda&= u_{m_a}\\
  &=u_{m_d}+ (u_{m_b}- u_{m_c})+ 2\cdot (u_{m_b}- u_{m_d})\\
  &=u_{m_d}+ 2\cdot (u_{m_b}- u_{m_c})+ (u_{m_b}- u_{m_d})\\
  &=u_{m_b}+ (u_{m_b}- u_{m_c})+ (u_{m_b}- u_{m_d}).
\end{align*}

If $M$ denotes the multiple of the $x$-vector that we add, then $M$
has to satisfy the equality
\begin{align*}
  \mu&=u_{m_f}+(u_{m_b}- u_{m_c})\cdot m_{fy}+ (u_{m_b}-
  u_{m_d})\cdot m_{fz}+3\cdot M\\
  &=u_{m_e}+(u_{m_b}- u_{m_c})\cdot m_{ey}+ (u_{m_b}-
  u_{m_d})\cdot m_{ez}+M
\end{align*}
where $m_{fy}$ is the second coordinate of $m_f$ etc., so that then
the new coefficients of $m_e$ and $m_f$ are also equal.
So we can solve for $M$ and then express the distance
$\frac{\lambda-\mu}{3}$ of the singular point from $V$ as
\begin{equation}
         \label{eq:b.1.1.4}
         \begin{aligned}
        \frac{\lambda-\mu}{3}=&
  \frac{u_{m_a}}{3}-
  \left(\frac{u_{m_e}}{2}-\frac{u_{m_f}}{6}\right)-(u_{m_b}-u_{m_c})\cdot
  \left(\frac{m_{ey}}{2}-\frac{m_{fy}}{6}\right) \\
  &-(u_{m_b}-u_{m_d})\cdot \left(\frac{m_{ez}}{2}-\frac{m_{fz}}{6}\right).
\end{aligned}
\end{equation}

This settles case (b.1) of Theorem \ref{thm2}.

%
%
%

\subsection{Weight class as in Lemma \ref{lem-chains}(b), circuit (D) of Remark \ref{rem-circuit}}
Let $F_{s-4}'=\{a,b,c,d\}$, $F_j'=\{e,f\}$, and assume first that
$m_e$ and $m_f$ lie on
different sides of the plane spanned by
$m_a,m_b,m_c$ and $m_d$. Since the two points  $m_e$ and $m_f$ have
the biggest heights of points outside the plane, it follows that
both form a pyramid with $m_a,m_b,m_c$ and $m_d$ in the subdivision.
By assumption both pyramids cannot have any lattice point besides
the five vertices. It follows from Lemma 3.3 of \cite{Sca85} that
the lattice distance of both points to the plane is one. Now solve
the equations given by the tropical polytope to get the positions of
the two vertices corresponding to the two pyramids. Without
restriction, we can assume that $m_a,m_b,m_c$ and $m_d$ lie in the
$x=0$-plane, it follows that the $x$-coordinate of $m_e$ is $-1$ and
the $x$-coordinate of $m_f$ is $1$. Since $m_a,m_b,m_c$ and $m_d$
have biggest and equal height, it follows that the edge dual to the
convex hull of $m_a,m_b,m_c$ and $m_d$ satisfies the equations $y=0$
and $z=0$. If $\lambda=u_{m_a}$ is the biggest weight (the weight of
$m_a,m_b,m_c$ and $m_d$), and $\mu=u_{m_e}$ is the weight of $m_e$
and $m_f$, it follows that the vertex dual to the pyramid with
vertex $m_e$ is at $(\mu-\lambda,0,0)$ and the vertex dual to the
pyramid with vertex $m_f$ is at $(\lambda-\mu,0,0)$. The singular
point $(0,0,0)$ is thus exactly in the middle of the edge dual to
the convex hull of $m_a,m_b,m_c$ and $m_d$.

Now assume $m_e$ and $m_f$ lie on the same side of the plane spanned
by $m_a,m_b,m_c$ and $m_d$. It follows from Lemma 3.3 of \cite{Sca85}
again that none of these two points can have an integral distance
larger than one to the plane, or it would form a pyramid with interior
lattice points. Thus both $m_e$ and $m_f$ have integral distance one,
and form a ``\emph{triangular roof}'' with $m_a,m_b,m_c$ and $m_d$.
Again, then the dual subdivision does not correspond to a cone of codimension $1$ of the secondary fan, and we do not consider the situation.
pyramid with base $m_a,m_b,m_c$ and $m_d$. We can again solve for
the position of vertex dual to this pyramid and get
$(\nu-\lambda,0,0)$. The singular point is on an edge ending at a
vertex $V_1$ adjacent to $5$ edges and $9$ $2$-dimensional polyhedra
and at a vertex $V_2$ with $5$ adjacent edges and $8$
$2$-dimensional polyhedra, and its distance to $V_2$ is bigger or
equal to its distance to $V_1$. This settles case (b.2) of Theorem
\ref{thm2}.

\subsection{Weight class as in Lemma \ref{lem-chains}(c), circuit (E) of Remark \ref{rem-circuit}}
\label{sec-barycenter}
With the notation from Lemma \ref{lem-chains}(c) let
$F_{s-4}'=\{a,b,c\}$ and $F_j'=\{d,e,f\}$. We may assume that
$m_a=(0,0,0)$, $m_b=(0,0,1)$, and $m_c=(0,0,2)$. We then
distinguish two cases. Either there is no plane containing the
$z$-axis such that $m_d$, $m_e$ and $m_f$ are all on one side of the
plane, or there is such a plane.

\subsubsection{Assume there is no plane through the $z$-axis with
  $m_d$, $m_e$, and $m_f$ all on the same side of the plane.}\label{subsec:E1}

In a first step we want to classify the possible polytopes spanned by
$m_a,\ldots,m_f$, and then we will see how the corresponding tropical
surfaces look like locally at the singular point.

\begin{lemma}\label{lem:nointlatpt}
  Let $P=\conv\big((0,0,0),(0,0,2),m,m'\big)$ with $m,m'\in\Z^3$ be
  a $3$-dimensional lattice polytope such that
  \begin{equation}
    \label{eq:nointlatpt:1}
    P\cap\Z^3=\{(0,0,0),(0,0,1),(0,0,2),m,m'\}.
  \end{equation}
  Projecting $P$ orthogonally onto the $xy$-plane
  we get a triangle $T$ which contains no interior lattice point and
  where the edges with vertex $(0,0)$ contain no relative interior point.
\end{lemma}
\begin{proof}
  We denote by $\pi:P\longrightarrow\R^2:(x,y,z)\mapsto (x,y)$
  the orthogonal projection onto the $xy$-plane, so that
  $T=\pi(P)$.

  Applying a suitable coordinate change in $\Gl_3(\Z)$ we may assume that
  $m'=(0,\beta',\gamma')$ and $m=(\alpha,\beta,\gamma)$ with $\beta'>0$. If $\beta'>1$
  then $\pi^{-1}(0,1)$ is
  a line segment of Euclidean length at least one and  it thus contains a
  lattice point in contradiction to
  \eqref{eq:nointlatpt:1}. Applying a coordinate change again we can assume $0\leq \beta<\alpha$. Since $\beta'=1$ the
  edge of $T$ connecting the vertex $(0,0)$ with $(0,\beta')$ has no relative interior point. If $\beta=0$ or $\beta=1$ the
  statement holds obviously, since then $T$ is a triangle of lattice
  height one (see Figure \ref{fig:trianglesheightone}).
  Note here that for $\beta=0$ necessarily $\alpha=1$ since otherwise
  above $\pi^{-1}(1,0)$ would contain an interior lattice point.
  \begin{figure}[h]
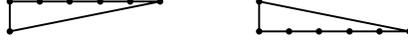

    \centering
    \begin{texdraw}
      \drawdim cm \relunitscale 0.4
      \move (0 0) \lvec (0 1) \lvec (5 1) \lvec (0 0)
      \move (0 0) \fcir f:0 r:0.1
      \move (0 1) \fcir f:0 r:0.1
      \move (1 1) \fcir f:0 r:0.1
      \move (2 1) \fcir f:0 r:0.1
      \move (3 1) \fcir f:0 r:0.1
      \move (4 1) \fcir f:0 r:0.1
      \move (5 1) \fcir f:0 r:0.1
    \end{texdraw}
    \hspace*{1cm}
    \begin{texdraw}
      \drawdim cm \relunitscale 0.4
      \move (0 0) \lvec (0 1) \lvec (5 0) \lvec (0 0)
      \move (0 0) \fcir f:0 r:0.1
      \move (0 1) \fcir f:0 r:0.1
      \move (1 0) \fcir f:0 r:0.1
      \move (2 0) \fcir f:0 r:0.1
      \move (3 0) \fcir f:0 r:0.1
      \move (4 0) \fcir f:0 r:0.1
      \move (5 0) \fcir f:0 r:0.1
    \end{texdraw}
    \caption{Lattice triangles of lattice height one.}
    \label{fig:trianglesheightone}
  \end{figure}

  We may therefore assume
  \begin{equation}
    \label{eq:nointlatpt:2}
    m'=(0,1,\gamma') \text{ and } m=(\alpha,\beta,\gamma)
    \text{ with } 2\leq \beta<\alpha.
  \end{equation}
  Moreover, we must have $\gcd(\alpha,\beta)=1$, since $\alpha=k\cdot d$ and
  $\beta=l\cdot d$ with $d\geq 2$ would
  imply that $\pi^{-1}(k,l)$ is a line segment of lattice length at
  least one and thus contains a lattice point in contradiction to
  \eqref{eq:nointlatpt:1}, see Figure \ref{fig:fibrewithlatticepoint}.
  \begin{figure}[h]
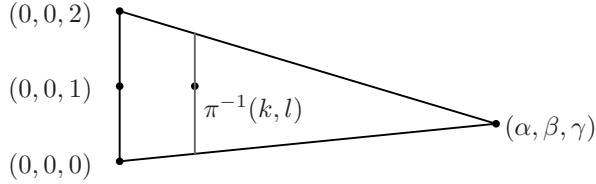

    \centering
    \begin{texdraw}
      \drawdim cm \relunitscale 0.5
      \move (0 0) \lvec (0 4) \lvec (10 1) \lvec (0 0)
      \move (0 0) \fcir f:0 r:0.1
      \move (0 2) \fcir f:0 r:0.1
      \move (0 4) \fcir f:0 r:0.1
      \move (10 1) \fcir f:0 r:0.1
      \move (2 2) \fcir f:0 r:0.1
      \htext (-3 -0.5){$(0,0,0)$}
      \htext (-3 1.5){$(0,0,1)$}
      \htext (-3 3.5){$(0,0,2)$}
      \htext (10.2 0.5){$(\alpha,\beta,\gamma)$}
      \htext (2.2 1){$\pi^{-1}(k,l)$}
      \setgray 0.4
      \move (2 0.2) \rlvec (0 3.2)
    \end{texdraw}
    \caption{$\pi^{-1}(k,l)$ contains a lattice point.}
    \label{fig:fibrewithlatticepoint}
  \end{figure}
  Therefore, also the edge of $T$ connecting vertex $(0,0)$ with $(\alpha,\beta)$ has no
  relative interior point, and if we divide $\alpha$ by $\beta$ with
  remainder we get
  \begin{equation}
    \label{eq:nointlatpt:3}
    \alpha=q\cdot \beta+r \text{ with } 1\leq r\leq \beta-1 \text{ and } q\geq 1.
  \end{equation}
  The triangle $T$ can be described by inequalities as follows
  \begin{displaymath}
    T=\left\{(x,y)\;\Big|\;x\geq 0,\; y\geq\frac{\beta}{\alpha}\cdot x,\;
      y\leq\frac{\beta-1}{\alpha}\cdot x+1\right\},
  \end{displaymath}
  which ensures that
  \begin{displaymath}
    (q,1)\in T.
  \end{displaymath}
  We now want to show that
  \begin{displaymath}
    \pi^{-1}(q,1)\cap\Z^3\not=\emptyset,
  \end{displaymath}
  which will be a contradiction to \eqref{eq:nointlatpt:1}.

  An easy computation shows that
  \begin{displaymath}
    \pi^{-1}(q,1)=\left\{\left(q,1,\frac{q\cdot \gamma+r\cdot \gamma'+z}{q\cdot \beta+r}\right)
      \;\Big|\; 0\leq z\leq 2\cdot q\cdot(\beta-1)\right\},
  \end{displaymath}
  and we have to show that there is a $0\leq z\leq 2\cdot q\cdot (\beta-1)$ such that
  \begin{equation}
    \label{eq:nointlatpt:4}
    q\cdot \beta+r\;\big|\; (q\cdot \gamma+r\cdot \gamma')+z.
  \end{equation}
  We consider first the special case $\beta=2$. Then necessarily $r=1$ and
  there is of course a $0\leq z\leq 2\cdot q$ such that $q\cdot
  \beta+r=2\cdot q+1$ divides $(q\cdot \gamma+\gamma')+z$.

  Next we consider the special case $(q,r)=(1,\beta-1)$, and we have to
  check if $q\cdot \beta+r=2\cdot \beta-1$ divides $(\gamma+(\beta-1)\cdot \gamma')+z$ for some $0\leq
  z\leq 2\cdot \beta-2$, which is obviously the case.

  For the general case we may now assume that $\beta\geq 3$ and
  $(q,r)\not=(1,\beta-1)$. Taking \eqref{eq:nointlatpt:2} and
  \eqref{eq:nointlatpt:3} into account it follows that
  \begin{displaymath}
    \beta\geq 2+\frac{r}{q},
  \end{displaymath}
  or equivalently
  \begin{displaymath}
    2\cdot q\cdot (\beta-1)\geq q\cdot \beta+r.
  \end{displaymath}
  But then, there is definitely a $0\leq z\leq 2\cdot q\cdot (\beta-1)$
  such that \eqref{eq:nointlatpt:4} is satisfied.

  So the case $2\leq\beta<\alpha$ cannot occur, and this finishes the
  proof.
\end{proof}

\begin{proposition}\label{prop:classification(E)1}
  Let $P$ be a lattice polytope which is the
  convex hull of a circuit of type (E) and three additional lattice
  points $m$, $m'$ and $m''$ such that any two of these together with
  the circuit span $\R^3$,
  $P$ contains only the given six lattice points, and there is no
  plane through the $z$-axis such that $m$, $m'$ and $m''$ are all on
  the same side of the plane, see Figure \ref{fig:polytope(E)1}.
  \begin{figure}[h]
    \centering
    \includegraphics[width=8cm]{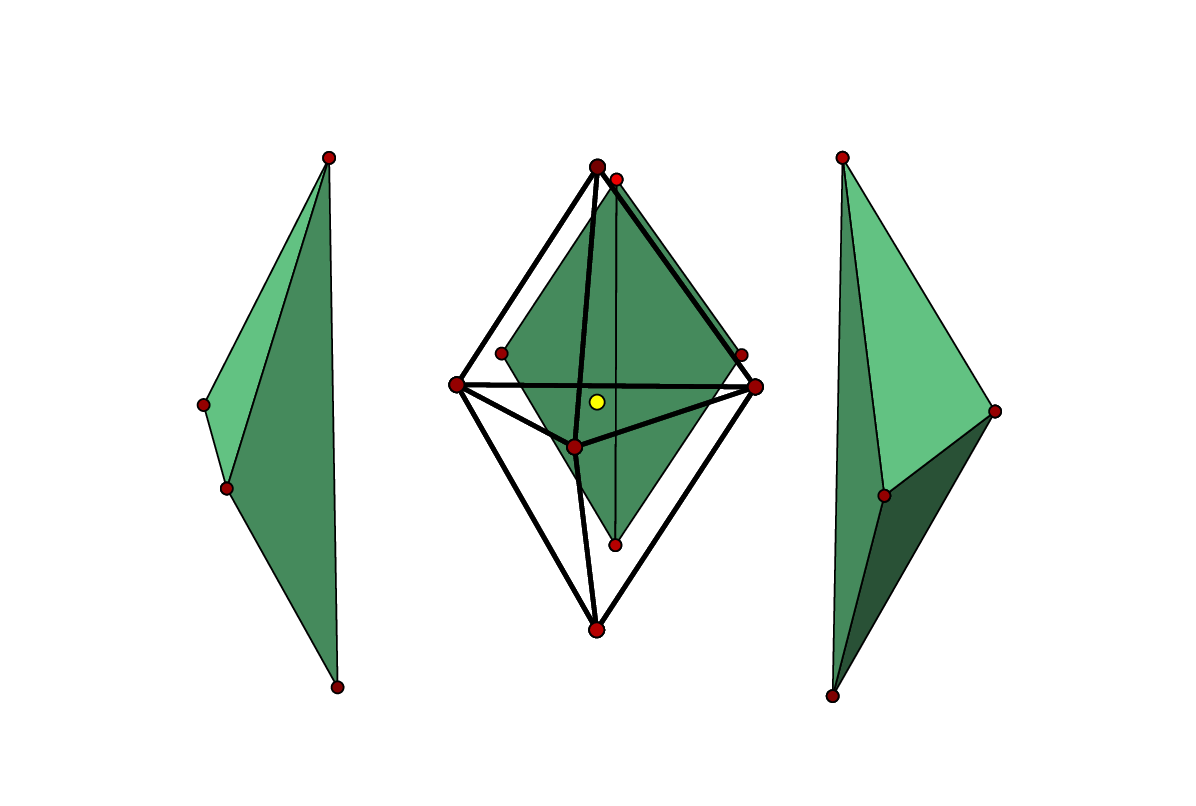}
    \caption{A lattice polytope $P$ as in Proposition
      \ref{prop:classification(E)1} with subdivision.}
    \label{fig:polytope(E)1}
  \end{figure}

  Then 
  the circuit is given up to IUA-equivalence
by $(0,0,0)$, $(0,0,1)$, and $(0,0,2)$, and the lattice points $m$,
  $m'$, and $m''$ satisfy the conditions in exactly one of the following
  cases:
  \begin{enumerate}
  \item $m=(0,1,\gamma)$, $m'=(1,0,\gamma')$, and $m''=(-1,-1,\gamma'')$ with
    $\gamma,\gamma',\gamma''\in\Z$ arbitrary.
  \item $m=(0,1,\gamma)$, $m'=(2,1,\gamma')$, and $m''=(-1,-1,\gamma'')$ with
    $\gamma,\gamma',\gamma''\in\Z$ such that $\gamma\not\equiv \gamma'\;(\mod 2)$.
  \item $m=(0,1,\gamma)$, $m'=(3,1,\gamma')$, and $m''=(-1,-1,\gamma'')$ with
    $\gamma,\gamma',\gamma''\in\Z$ such that $\gamma\not\equiv \gamma'\;(\mod 3)$ and $\gamma'\not\equiv
    \gamma''\;(\mod 2)$.
  \item $m=(0,1,\gamma)$, $m'=(3,1,\gamma')$, and $m''=(-3,-2,\gamma'')$ with
    $\gamma,\gamma',\gamma''\in\Z$ such that $\gamma\not\equiv \gamma'\not\equiv \gamma''\not\equiv \gamma\;(\mod
    3)$.
  \end{enumerate}
  \begin{center}
  \end{center}
\end{proposition}
\begin{proof}
 It is clear that 
  the circuit (E) 
  is IUA-equivalent to $(0,0,0)$, $(0,0,1)$, and $(0,0,2)$.
  If we denote by
  $\pi:P\longrightarrow\R^2:(x,y,z)\mapsto(x,y)$ the projection
  onto the $xy$-plane then
  $\pi(P)$ is a triangle which decomposes into three triangles
  $\pi(P)=T\cup T'\cup T''$ as
  in Lemma \ref{lem:nointlatpt}, see Figure~\ref{fig:TT'T''}.
  \begin{figure}[h]
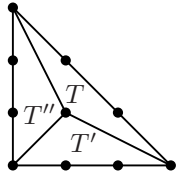

    \centering
    \begin{texdraw}
      \drawdim cm \relunitscale 0.7
      \move (0 0) \lvec (3 0) \lvec (0 3) \lvec (0 0)
      \move (0 0) \fcir f:0 r:0.1
      \move (1 0) \fcir f:0 r:0.1
      \move (2 0) \fcir f:0 r:0.1
      \move (3 0) \fcir f:0 r:0.1
      \move (0 1) \fcir f:0 r:0.1
      \move (1 1) \fcir f:0 r:0.1
      \move (2 1) \fcir f:0 r:0.1
      \move (0 2) \fcir f:0 r:0.1
      \move (1 2) \fcir f:0 r:0.1
      \move (0 3) \fcir f:0 r:0.1
      \move (0 0) \lvec (1 1)
      \move (3 0) \lvec (1 1)
      \move (0 3) \lvec (1 1)
      \htext (1.1 0.25){$T'$}
      \htext (0.2 0.8){$T''$}
      \htext (1 1.2){$T$}
    \end{texdraw}
    \caption{$\pi(P)=T\cup T'\cup T''$ decomposes as a union of three triangles.}
    \label{fig:TT'T''}
  \end{figure}
  Lemma \ref{lem:nointlatpt} therefore implies that $(0,0)$ is the only
  interior lattice point of $\pi(P)$.
  Lattice polygons with exactly one interior
  lattice point have been classified up to IUA-equivalence, see e.g.\ \cite{Rab89} or \cite{PR00}, and among
  them are exactly five triangles as shown in Figure
  \ref{fig:latticetrianglesoneintpt}, where the interior lattice point
  is $(0,0)$.
  \begin{figure}[h]
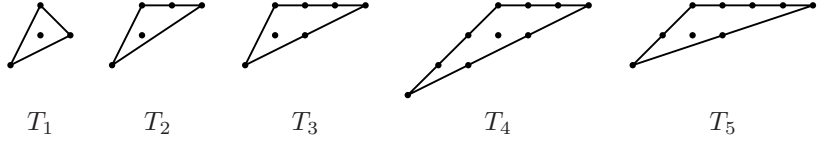

    \centering
     \begin{tabular}[m]{c@{\hspace*{0.5cm}}c@{\hspace*{0.5cm}}c@{\hspace*{0.5cm}}c@{\hspace*{0.5cm}}c}
       \begin{texdraw}
         \drawdim cm \relunitscale 0.4
         \move (0 1) \lvec (1 0) \lvec (-1 -1) \lvec (0 1)
         \move (0 0) \fcir f:0 r:0.1
         \move (0 1) \fcir f:0 r:0.1
         \move (1 0) \fcir f:0 r:0.1
         \move (-1 -1) \fcir f:0 r:0.1
         \move (0 -2) \fcir f:1 r:0.1
       \end{texdraw}
       &
       \begin{texdraw}
         \drawdim cm \relunitscale 0.4
         \move (0 1) \lvec (2 1) \lvec (-1 -1) \lvec (0 1)
         \move (0 0) \fcir f:0 r:0.1
         \move (0 1) \fcir f:0 r:0.1
         \move (1 1) \fcir f:0 r:0.1
         \move (2 1) \fcir f:0 r:0.1
         \move (-1 -1) \fcir f:0 r:0.1
         \move (0 -2) \fcir f:1 r:0.1
       \end{texdraw}
       &
       \begin{texdraw}
         \drawdim cm \relunitscale 0.4
         \move (0 1) \lvec (3 1) \lvec (-1 -1) \lvec (0 1)
         \move (0 0) \fcir f:0 r:0.1
         \move (0 1) \fcir f:0 r:0.1
         \move (1 1) \fcir f:0 r:0.1
         \move (2 1) \fcir f:0 r:0.1
         \move (3 1) \fcir f:0 r:0.1
         \move (1 0) \fcir f:0 r:0.1
         \move (-1 -1) \fcir f:0 r:0.1
         \move (0 -2) \fcir f:1 r:0.1
       \end{texdraw}
       &
       \begin{texdraw}
         \drawdim cm \relunitscale 0.4
         \move (0 1) \lvec (3 1) \lvec (-3 -2) \lvec (0 1)
         \move (0 0) \fcir f:0 r:0.1
         \move (0 1) \fcir f:0 r:0.1
         \move (1 1) \fcir f:0 r:0.1
         \move (2 1) \fcir f:0 r:0.1
         \move (3 1) \fcir f:0 r:0.1
         \move (1 0) \fcir f:0 r:0.1
         \move (-1 -1) \fcir f:0 r:0.1
         \move (-3 -2) \fcir f:0 r:0.1
         \move (-2 -1) \fcir f:0 r:0.1
         \move (-1 0) \fcir f:0 r:0.1
       \end{texdraw}
       &
       \begin{texdraw}
         \drawdim cm \relunitscale 0.4
         \move (0 1) \lvec (4 1) \lvec (-2 -1) \lvec (0 1)
         \move (0 0) \fcir f:0 r:0.1
         \move (0 1) \fcir f:0 r:0.1
         \move (1 1) \fcir f:0 r:0.1
         \move (2 1) \fcir f:0 r:0.1
         \move (3 1) \fcir f:0 r:0.1
         \move (4 1) \fcir f:0 r:0.1
         \move (1 0) \fcir f:0 r:0.1
         \move (-2 -1) \fcir f:0 r:0.1
         \move (-1 0) \fcir f:0 r:0.1
         \move (0 -2) \fcir f:1 r:0.1
       \end{texdraw}
       \\
       $T_1$&$T_2$&$T_3$&$T_4$&$T_5$ \\[0.4cm]
     \end{tabular}
    \caption{The five lattice triangles with one interior lattice point.}
    \label{fig:latticetrianglesoneintpt}
  \end{figure}
  Applying a $\Z$-linear coordinate change we may therefore assume
  that $\pi(P)$ is one of these five triangles. In each of
  the cases it remains to check whether there exist
  polytopes $P$ that project to the triangle and what
  restrictions this poses on the third component of the lattice points
  $m$, $m'$, and $m''$. Actually, the only obstruction is that above
  the relative interior lattice points on the edges of the triangles
  there should be no lattice point in $P$. If such an edge has
  $k$ relative interior lattice points and the $z$-coordinates of the
  vertices of the edge differ by $l$, then some of the relative
  interior lattice points lifts to a lattice point if and only if
  $k+1$ and $l$  are not coprime.
  Therefore, $T_1,\ldots,T_4$ lead to the
  four cases mentioned in the statement of the proposition.
  For $T_5$ we would need points $m=(0,1,\gamma)$, $m'=(4,1,\gamma')$, and
  $m''=(-2,-1,\gamma'')$ such that each of the differences $\gamma-\gamma'$, $\gamma-\gamma''$ and $\gamma'-\gamma''$
  is coprime to two. That is obviously not possible, so
  that $T_5$ cannot be the projection of any $P$.
\end{proof}

\begin{figure}[h]
  \centering
  \includegraphics[width=5cm]{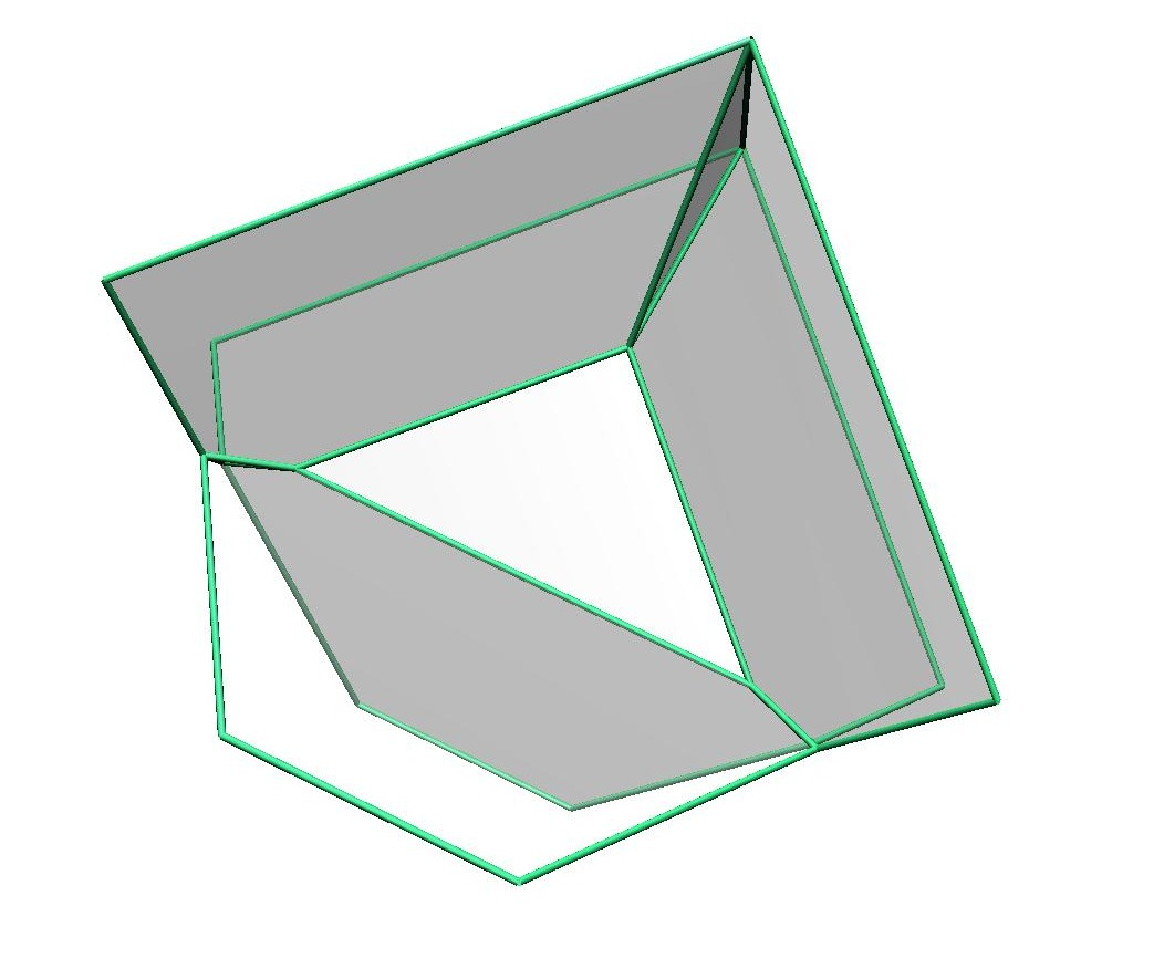}
  \caption{The triangle in the tropical surface dual to the circuit.}
  \label{fig:tropsurftriangle}
\end{figure}

\begin{figure}[h]
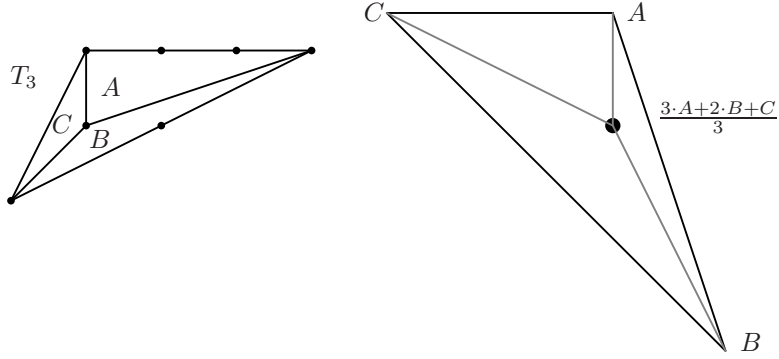

  \centering
  \begin{texdraw}
    \drawdim cm \relunitscale 1
    \move (0 1) \lvec (3 1) \lvec (-1 -1) \lvec (0 1)
    \move (0 0) \fcir f:0 r:0.05
    \move (0 1) \fcir f:0 r:0.05
    \move (1 1) \fcir f:0 r:0.05
    \move (2 1) \fcir f:0 r:0.05
    \move (3 1) \fcir f:0 r:0.05
    \move (1 0) \fcir f:0 r:0.05
    \move (-1 -1) \fcir f:0 r:0.05
    \move (0 -2) \fcir f:1 r:0.05
    \move (0 0) \lvec (-1 -1)
    \move (0 0) \lvec (0 1)
    \move (0 0) \lvec (3 1)
    \htext (0.2 0.4){$A$}
    \htext (0.05 -0.3){$B$}
    \htext (-0.45 -0.1){$C$}
    \htext (-1 0.5){$T_3$}

    \move (7 1.5) \lvec (8.5 -3) \lvec (4 1.5) \lvec (7 1.5)
    \htext (7.2 1.4){$A$}
    \htext (8.7 -3){$B$}
    \htext (3.7 1.4){$C$}
    \htext (7.6 -0.1){$\frac{3\cdot A+2\cdot B+C}{3}$}
    \move (7 0) \fcir f:0 r:0.1
    \setgray 0.5
    \lvec (7 1.5)
    \move (7 0) \lvec (8.5 -3)
    \move (7 0) \lvec (4 1.5)
  \end{texdraw}
  \caption{$T_3=\pi(\Delta)$ and the dual triangle in the tropical
    surface showing $(0,0,0)$ as the weighted barycenter $\frac{3\cdot
    A+2\cdot B+C}{3}$.}
  \label{fig:baricenter}
\end{figure}

In order to understand how the tropicalisation of the singular point
locally looks like in the case we are considering, assume first that
the subdivision contains a polytope as considered in
Proposition~\ref{prop:classification(E)1}, and it is subdivided into
the three polytopes $\Delta_1=\conv(m_a,m_c,m_d,m_e)$,
$\Delta_2=\conv(m_a,m_c,m_d,m_f)$ and
$\Delta_3=\conv(m_a,m_c,m_e,m_f)$, see Figure \ref{fig:polytope(E)1}.
The circuit
$\{m_a,m_b,m_c\}$ is then dual to a triangle in the tropical
surface whose vertices are dual to $\Delta_1$, $\Delta_2$, and
$\Delta_3$, see Figure \ref{fig:tropsurftriangle}.
We assume as before that $m_a=(0,0,0)$, $m_b=(0,0,1)$ and
$m_c=(0,0,2)$. Recall that we can project $P$ to the $(x,y)$-plane
and obtain three triangles $T$, $T'$ and $T''$ as in Figure
\ref{fig:TT'T''}. The midpoint is $(0,0)$. Denote the coordinates of
the three vertices by $(r_1,s_1)$, $(r_2,s_2)$ and $(r_3,s_3)$. Let
us use the tropical polynomial to solve for the coordinates
$(x,y,z)$ of the three vertices dual to $\Delta_1$, $\Delta_2$, and
$\Delta_3$. By assumption the heights associated to the lattice
points satisfy $u_{m_a}=u_{m_b}=u_{m_c}$ and
$u_{m_d}=u_{m_e}=u_{m_f}$, and we set $u=u_{m_a}-u_{m_d}$. For any
$i=1,2,3$, the equation $u+z=u$ has to be satisfied, so any of the
three vertices has $z$-coordinate $0$. In fact, the whole triangle
dual to the circuit satisfies $z=0$. So we only have to solve for
the $(x,y)$-coordinates of the vertices. For any choice of
$(i,j)=(1,2), (2,3)$ or $(3,1)$, the vertex dual to the polytope
which projects  to the triangle spanned by $(0,0)$, $(r_i,s_i)$
and
$(r_j,s_j)$ has to satisfy the equations $u=r_ix+s_iy$ and
$u=r_jx+s_jy$, which are solved by
$(x,y)=\frac{1}{r_is_j-s_ir_j}\cdot (s_ju-s_iu,r_iu-r_ju)$. Now
assign to each of the vertices the area of the projection of the
dual polytope, i.e.\ $(r_is_j-s_ir_j)$, as weight. Then it follows
that the weighted sum of the three vertices is $(0,0,0)$, i.e.\ the
singular point. Thus, the singular point tropicalises precisely to
the \emph{weighted barycenter} of the triangle dual to the circuit.
Figure \ref{fig:baricenter} depicts this situation for the case that
the projection is the triangle $T_3$ of Figure
\ref{fig:latticetrianglesoneintpt}.

\begin{figure}[h]
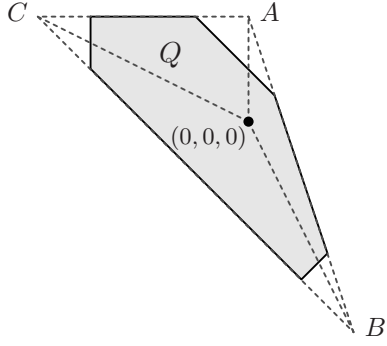

  \centering
  \begin{texdraw}
    \drawdim cm \relunitscale 0.7
    \setgray 1
    \move (-3 2) \lvec (-1 2) \lvec (0.5 0.5) \lvec (1.5 -2.5)
    \lvec (1 -3) \lvec (-3 1) \lvec (-3 2) \lfill f:0.9
    \setgray 0
    \move (-3 2) \lvec (-1 2) \lvec (0.5 0.5) \lvec (1.5 -2.5)
    \lvec (1 -3) \lvec (-3 1) \lvec (-3 2)
    \setgray 0.3
    \lpatt (0.067 0.1)
    \move (0 2) \lvec (2 -4) \lvec (-4 2) \lvec (0 2)
    \htext (-1.7 1){\large $Q$}
    \htext (0.2 1.9){$A$}
    \htext (2.2 -4.1){$B$}
    \htext (-4.6 1.9){$C$}
    \htext (-1.5 -0.55){\small $(0,0,0)$}
    \move (0 0) \lvec (0 2)
    \move (0 0) \lvec (2 -4)
    \move (0 0) \lvec (-4 2)
    \move (0 0) \fcir f:0 r:0.1
  \end{texdraw}
  \caption{The origin as a generalised weighted barycenter.}
  \label{fig:baricenter2}
\end{figure}

If the subdivision locally around the circuit contains further lattice points, the local picture may look more complicated. However, the circuit
$\{m_a,m_b,m_c\}$ is still dual to a polygon $Q$ in the
$\{z=0\}$-plane. Moreover, in the subdivision there will still be polytopes
which contain $\conv(m_a,m_c,m_d)$ respectively $\conv(m_a,m_c,m_e)$
respectively $\conv(m_a,m_c,m_f)$ as a facet. Therefore, the
polygon $Q$ will have three edges dual to these facets. If one
computes the intersection points of the lines through these edges, one
gets three points $A$, $B$, and $C$ which would be dual to the
polytopes $\Delta_i$. This extension of the cell forms a virtual triangular cell, and the tropicalisation of the singular point is
still the weighted sum of the three points $A$, $B$ and $C$, see Figure
\ref{fig:baricenter2}. 

\begin{example}\label{ex:barycenter1}
  A concrete example for this behaviour is the singular point
  $H=(-1,-1,0)$ on the tropical surface in Example
  \ref{ex-thomas}. We have seen in Remark \ref{rem:ex-thomas} which
  weight class corresponds to the point $H=(-1,-1,0)$. We have
  \begin{displaymath}
    \begin{array}{lll}
      m_a=(0,0,0),& m_b=(0,0,1),& m_c=(0,0,2),  \\
      m_d=(-1,-1,0),&  m_e=(0,1,0), & m_f=(1,0,0),
    \end{array}
  \end{displaymath}
  and one further point $m_g=(1,2,1)$. The circuit $m_a,m_b,m_c$
  corresponds then to quadrangle $ABCD$ (see Figure \ref{fig:w-barycenter},
  \begin{figure}[h]
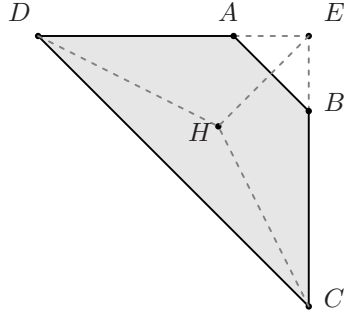

    \centering
    \begin{texdraw}
      \drawdim cm \relunitscale 0.2
      \move (5 -13) \lvec (-13 5) \lvec (0 5) \lvec (5 0)
      \lvec (5 -13) \lfill f:0.9
      \lpatt (0.3 0.5)
      \setgray 0.5
      \move (5 0) \lvec (5 5) \lvec (0 5)
      \move (5 0) \fcir f:0 r:0.2
      \move (0 5) \fcir f:0 r:0.2
      \move (5 5) \fcir f:0 r:0.2
      \move (-13 5) \fcir f:0 r:0.2
      \move (5 -13) \fcir f:0 r:0.2
      \move (-1 -1) \fcir f:0 r:0.2
      \htext (-15 6) {$D$}
      \htext (6 -13) {$C$}
      \htext (-1 6) {$A$}
      \htext (6 0) {$B$}
      \htext (6 6) {$E$}
      \htext (-3 -2) {$H$}
      \move (-13 5) \lvec (-1 -1) \lvec (5 -13)
      \move (5 5) \lvec (-1 -1)
    \end{texdraw}
    \caption{The singular point as barycenter.}
    \label{fig:w-barycenter}
  \end{figure}
  where the vertices $C=(5,-13,0)$ and $D=(-13,5,0)$ correspond to the
  polytopes $\Delta_C=\conv\{m_a,m_c,m_d,m_e\}$ respectively
  $\Delta_D=\conv\{m_a,m_c,m_d,m_f\}$ in the subdivision. The polytope
  $\Delta_E=\conv\{m_a,m_c,m_e,m_f\}$, however, is not part of the
  subdivision due to the presence of $m_g$ with an appropriate
  height. However, $\Delta_E$ defines a virtual point $E=(5,5,0)$, which is
  the intersection of the two lines determined by the facets
  $\conv\{m_a,m_c,m_e\}$ and $\conv\{m_a,m_c,m_f\}$ of $\Delta_C$
  respectively $\Delta_D$, and
  \begin{displaymath}
    H=\frac{1}{3}\cdot (C+D+E)
  \end{displaymath}
  is the barycenter of this virtual triangle in the tropical surface.
\end{example}

\subsubsection{Assume there is a plane through the $z$-axis with
  $m_d$, $m_e$, and $m_f$ all on the same side of the plane.}\label{subsubsec:virtualbarycenter}

Again we first want to classify the possible polytopes spanned by
$m_a,\ldots,m_f$, and then we will see how the corresponding tropical
surfaces look like locally at the singular point.

\begin{figure}[h]
  \centering
  \includegraphics[width=6cm]{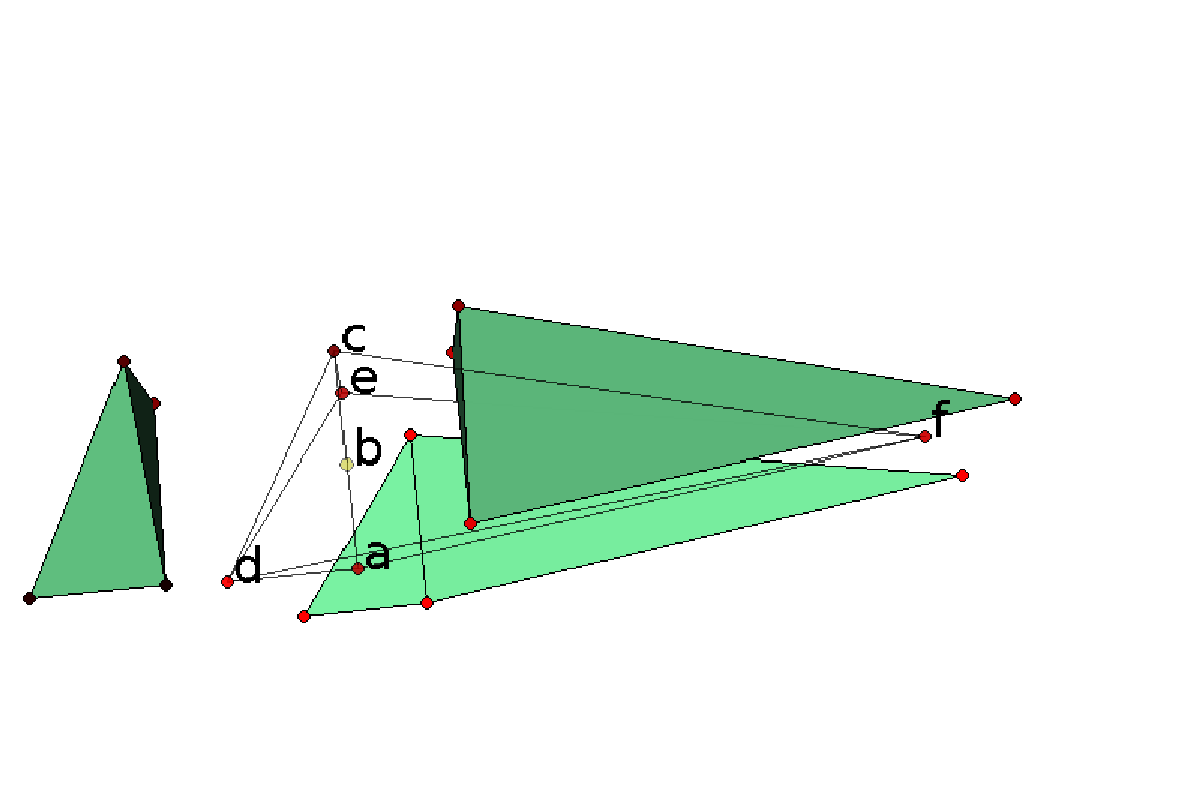}
  \includegraphics[width=6cm]{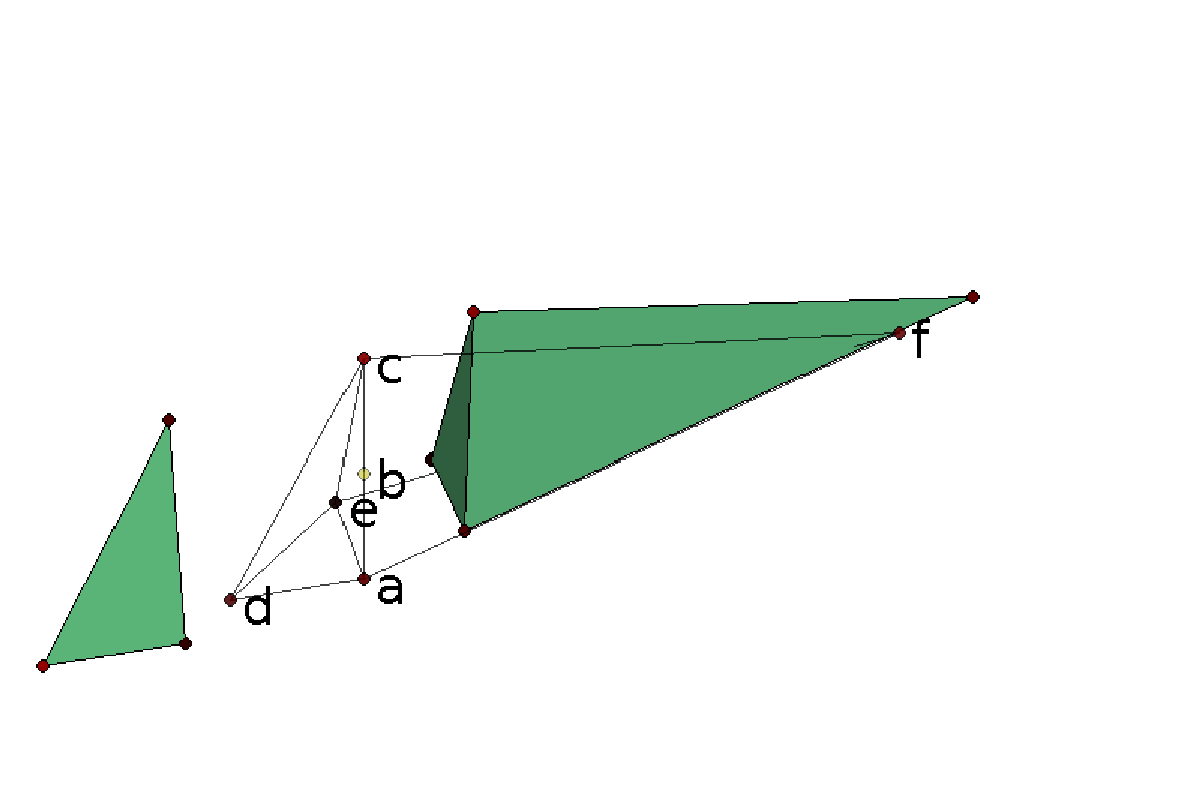}
  \caption{Possible subdivisions of $\Delta$.}
  \label{fig:subdivisionE2}
\end{figure}

\begin{proposition}\label{prop:classification(E)2}
  Let $P$ be a lattice polytope which is the
  convex hull of a circuit of type (E) and three additional lattice
  points $m$, $m'$, and $m''$ such that any two of these together with
  the circuit span $\R^3$,
  $P$ contains only the given six lattice points, and there is a
  plane through the $z$-axis such that $m$, $m'$ and $m''$ are all on
  the same side of the plane, see Figure \ref{fig:subdivisionE2}.

  Then 
the circuit is given up to IUA-equivalence
  by $(0,0,0)$, $(0,0,1)$, and $(0,0,2)$, and the lattice points $m$,
 $m'$, and $m''$ (up to reordering) satisfy the conditions in exactly one of the following
  cases:
  \begin{enumerate}
  \item $m=(-1,0,\gamma)$, $m'=(0,1,\gamma')$, and $m''=(\alpha'',1,\gamma'')$ with
    $\alpha''\geq 1$, $\gamma\in\Z$ arbitrary and $\gcd(\gamma''-\gamma',\alpha'')=1$.
  \item $m=(\alpha,1,\gamma)$, $m'=(\alpha+l,1,\gamma+k)$, and
    $m''=(\alpha+2\cdot l,1,\gamma+2\cdot k)$ with
    $\alpha,\gamma\in\Z$ arbitrary and $\gcd(l,k)=1$.
  \item $m=(\alpha,1,\gamma)$, $m'=(\alpha',1,\gamma')$, and $m''=(\alpha'',1,\gamma'')$ with
    $$\det\begin{pmatrix}\alpha'-\alpha&\alpha''-\alpha\\\gamma'-\gamma&\gamma''-\gamma\end{pmatrix}=\pm 1.$$
  \end{enumerate}
\end{proposition}
\begin{proof}
   Up to IUA-equivalence we may assume that the
  circuit is $(0,0,0)$, $(0,0,1)$, and $(0,0,2)$. Projecting $\Delta$
  to the $xy$-plane the points $\pi(m)$, $\pi(m')$, and $\pi(m'')$ lie
  in one half plane. Due to the assumptions on $\Delta$ no two of
  these points lie on the same line through the origin, and ordering
  these lines by their angle clockwise we may assume up to
  reordering that the
  points $\pi(m)$, $\pi(m')$, and $\pi(m'')$ come in this order, see
  Figure \ref{fig:clockwise} for possible configurations.
  \begin{figure}[h]
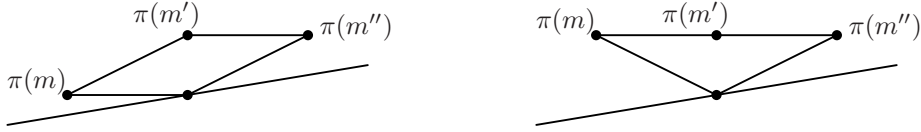

    \centering
    \bigskip
    \begin{texdraw}
      \drawdim cm \relunitscale 0.8
      \move (0 0) \fcir f:0 r:0.08
      \move (-3 -0.5) \lvec (3 0.5)
      \move (0 0) \lvec (-2 0) \lvec (0 1) \lvec (2 1) \lvec (0 0)
      \move (-2 0) \fcir f:0 r:0.08
      \move (0 1) \fcir f:0 r:0.08
      \move (2 1) \fcir f:0 r:0.08
      \htext (-3 0){$\pi(m)$}
      \htext (-0.9 1.1){$\pi(m')$}
      \htext (2.2 0.9){$\pi(m'')$}
    \end{texdraw}
    \hspace{2cm}
    \begin{texdraw}
      \drawdim cm \relunitscale 0.8
      \move (0 0) \fcir f:0 r:0.08
      \move (-3 -0.5) \lvec (3 0.5)
      \move (0 0) \lvec (-2 1) \lvec (0 1) \lvec (2 1) \lvec (0 0)
      \move (-2 1) \fcir f:0 r:0.08
      \move (0 1) \fcir f:0 r:0.08
      \move (2 1) \fcir f:0 r:0.08
      \htext (-3 1){$\pi(m)$}
      \htext (-0.9 1.1){$\pi(m')$}
      \htext (2.2 0.9){$\pi(m'')$}
    \end{texdraw}
    \caption{$\pi(\Delta)$ with the separating hyperplane.}
    \label{fig:clockwise}
  \end{figure}

  We should note here first that in $\pi(\Delta)$ the point $\pi(m')$
  cannot be an interior point of
  $\conv\big((0,0),\pi(m),\pi(m'')\big)$, since otherwise $\conv\big((0,0,0),(0,0,2),m,m''\big)$ will be in the
  subdivision of $\Delta$ which therefore satisfies the assumptions on
  Lemma \ref{lem:nointlatpt}, but $\pi(m')$ would violate these
  assumptions. It is then natural to distinguish the two cases that
  either $\pi(m')$ is on the line segment connecting $\pi(m)$ and
  $\pi(m'')$, i.e.\ $\pi(\Delta)$ is a triangle as shown on the right
  hand side of Figure \ref{fig:clockwise}, or $\pi(\Delta)$ is a
  quadrangle as shown on the left hand side of Figure
  \ref{fig:clockwise}. In any case, applying Lemma
  \ref{lem:nointlatpt} to the convex hull of the circuit and two of
  the further lattice points $m$, $m'$, and $m''$, we see that each of
  the points $\pi(m)$, $\pi(m')$, and $\pi(m'')$ has lattice distance
  one from the origin.

  Let us first consider the case that $\pi(\Delta)$ is a
  triangle. Up to IUA-equivalence we
  may assume that the line through $\pi(m)$, $\pi(m')$, and $\pi(m'')$
  is parallel to the $x$-axis, i.e.\ $\pi(m)=(\alpha,\beta)$,
  $\pi(m')=(\alpha',\beta)$, and $\pi(m'')=(\alpha'',\beta)$ with
  $\alpha<\alpha'<\alpha''$. By Lemma
  \ref{lem:nointlatpt} the triangle
  $\conv\big((0,0),\pi(m),\pi(m'')\big)$ has no interior lattice
  point and the number of lattice points on the boundary is
  $\alpha''-\alpha+2$, so that Pick's Formula implies $\beta=1$.
  \begin{figure}[h]
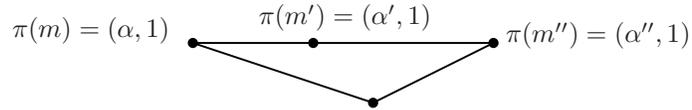

    \centering
    \begin{texdraw}
      \drawdim cm \relunitscale 0.8
      \move (0 0) \fcir f:0 r:0.08
      \move (0 0) \lvec (-3 1) \lvec (-1 1) \lvec (2 1) \lvec (0 0)
      \move (-3 1) \fcir f:0 r:0.08
      \move (-1 1) \fcir f:0 r:0.08
      \move (2 1) \fcir f:0 r:0.08
      \htext (-6 1){$\pi(m)=(\alpha,1)$}
      \htext (-1.9 1.2){$\pi(m')=(\alpha',1)$}
      \htext (2.2 0.9){$\pi(m'')=(\alpha'',1)$}
      \move (0 1.7)
      \fcir f:1 r:0.1
    \end{texdraw}
    \caption{The normal form of $\pi(\Delta)$ when it is a triangle.}
    \label{fig:casetriangle}
  \end{figure}

  This case now subdivides into two subcases, namely, that the points
  $m$, $m'$, and $m''$ lie on a line, respectively that they form a
  triangle. If the three points lie on a line,
  then $m'$ must be the
  midpoint of the line segment from $m$ to $m''$ and the line segment
  contains no further lattice point. Thus,
  $\gcd\big(\alpha''-\alpha,\gamma''-\gamma\big)=2$ is the only
  obstruction that has to be satisfied, and we are thus in Case (b) of
  the proposition with $l=\frac{\alpha''-\alpha}{2}$ and
  $k=\frac{\gamma''-\gamma}{2}$. If the three points $m$, $m'$, and
  $m''$ form a triangle, then the only obstruction to the condition
  that $\Delta$ contains no further lattice points is that this
  triangle should have lattice area one. This is precisely the
  condition of Case (c) in the proposition.

  It remains to consider the case that $\pi(\Delta)$ is a
  quadrangle.
As in the proof of Lemma \ref{lem:nointlatpt}, up to
  IUA-equivalence,
  $m'=(0,1,\gamma')$ and $m''=(\alpha'',\beta'',\gamma'')$ with $0\leq
  \beta''<\alpha''$. Moreover, since the triangle
  $T=\conv\big((0,0),\pi(m'),\pi(m'')\big)$ contains no interior lattice
  point due to Lemma \ref{lem:nointlatpt} Pick's Formula implies that
  $\beta''\in\{0,1\}$, and if $\beta''=1$ then necessarily
  $\alpha''=1$, since the lattice distance from $\pi(m'')$ to the
  origin is one. See Figure \ref{fig:quadranglecase}.
  \begin{figure}[h]
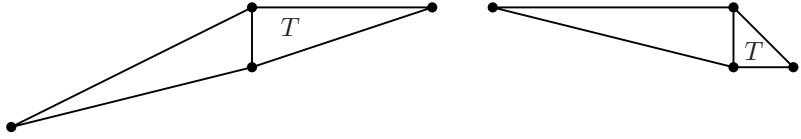

    \centering
    \begin{texdraw}
      \drawdim cm \relunitscale 0.8
      \move (0 0) \fcir f:0 r:0.08
      \move (0 0) \lvec (0 1) \lvec (3 1) \lvec (0 0) \lvec (-4 -1)
      \lvec (0 1)
      \move (3 1) \fcir f:0 r:0.08
      \move (0 1) \fcir f:0 r:0.08
      \move (-4 -1) \fcir f:0 r:0.08
      \htext (0.5 0.5){$T$}

      \move (8 0) \fcir f:0 r:0.08
      \move (8 0) \lvec (8 1) \lvec (9 0) \lvec (8 0) \lvec (4 1)
      \lvec (8 1)
      \move (9 0) \fcir f:0 r:0.08
      \move (8 1) \fcir f:0 r:0.08
      \move (4 1) \fcir f:0 r:0.08
      \htext (8.2 0.1){$T$}
    \end{texdraw}
    \caption{Possible configurations for the triangle $T=\conv\big((0,0),\pi(m'),\pi(m'')\big)$.}
    \label{fig:quadranglecase}
  \end{figure}

  Let us now consider the case $\beta''=1$ in more
  detail. The point $\pi(m)=(\alpha,\beta)$ has to lie below the line
  $\{y=1\}$ and above the line $\{\alpha\cdot y=x\}$. Thus $0\geq
  \beta>\alpha$, and applying Pick's Formula once again we find
  $\beta=0$, and then necessarily $\alpha=-1$. Analogously, we get in
  the case $\beta''=0$ that $\beta''=1$ and $\alpha\geq 1$. That is,
  $\pi(\Delta)$ is one of the quadrangles shown in Figure \ref{fig:quadranglecase2}.
  \begin{figure}[h]
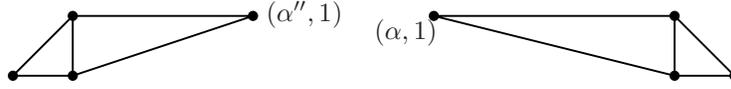

    \centering
    \begin{texdraw}
      \drawdim cm \relunitscale 0.8
      \move (0 0) \fcir f:0 r:0.08
      \move (0 0) \lvec (0 1) \lvec (3 1) \lvec (0 0) \lvec (-1 0)
      \lvec (0 1)
      \move (3 1) \fcir f:0 r:0.08
      \move (0 1) \fcir f:0 r:0.08
      \move (-1 0) \fcir f:0 r:0.08
      \htext (3.2 0.8){$(\alpha'',1)$}

      \move (10 0) \fcir f:0 r:0.08
      \move (10 0) \lvec (10 1) \lvec (11 0) \lvec (10 0) \lvec (6 1)
      \lvec (10 1)
      \move (11 0) \fcir f:0 r:0.08
      \move (10 1) \fcir f:0 r:0.08
      \move (6 1) \fcir f:0 r:0.08
      \htext (5 0.5){$(\alpha,1)$}
    \end{texdraw}
    \caption{The normal forms of $\pi(\Delta)$ when it is a quadrangle.}
    \label{fig:quadranglecase2}
  \end{figure}

  Obviously, reflecting at the plane $\{x=0\}$ and exchanging $m$ and
  $m''$ the two possible configuration types are equivalent, so that
  we may assume that $\beta=1$. We thus have $m=(-1,0,\gamma)$,
  $m'=(0,1,\gamma')$, and $m''=(\alpha'',1,\gamma'')$. Only above the
  line segment joining $\pi(m')$ and $\pi(m'')$ there could be an
  additional lattice point in $\Delta$ if the coordinates $\gamma'$
  and $\gamma''$ are chosen inappropriately, and the condition to avoid
  this is $\gcd(\gamma''-\gamma',\alpha'')=1$. We are thus in Case (a)
  of the proposition, and this finishes the proof.
\end{proof}

\longer{\begin{remark}\label{rem:defective}
  If $\Delta$ was the polytope $P$ in Proposition
  \ref{prop:classification(E)2} (b), then $\Delta$ would be defective
  (see e.g.\ Example 2.12 in \cite{DT10}), and we would not consider
  $\Delta$ at all. If $\Delta$, however, contains further lattice
  points besides those in $P$, then $\Delta$ need not be
  defective. But the weight class $C$ in $\Trop(\Ker(A))$
  corresponding to this situation will still be defective, as we will
  see further down when considering this case.
\end{remark}}

We now have to see how the tropical surface looks locally at the
tropicalisation of the singular point, i.e.\ locally at $(0,0,0)$.
As in Subsection \ref{subsec:E1} we want to restrict first to the case
where the Newton polytope $\Delta$ is just the convex hull of
$m_a,\ldots,m_f$, and in the notation of Proposition
\ref{prop:classification(E)2} we may assume that $m_d=m$, $m_e=m'$,
and $m_f=m''$.  Moreover, we will consider the Case (a) in Proposition
\ref{prop:classification(E)2} first.
In the subdivision of $\Delta$ there will be exactly
two polytopes which contain the circuit $m_a$, $m_b$, and $m_c$,
namely $\Delta_A=\conv(m_a,m_b,m_c,m_d,m_e)$ and
$\Delta_B=\conv(m_a,m_b,m_c,m_e,m_f)$, see Figure \ref{fig:subdivisionE2}. The
subdivision may contain a third polytope $\conv(m_a,m_d,m_e,m_f)$
respectively $\conv(m_c,m_d,m_e,m_f)$ which does not contain the circuit,
and which consequently will not matter for the singular point.

The tropicalisation of the singular point will then be contained in
the plane segment dual to the circuit. This segment will be unbounded,
but it has two vertices $A$ and $B$ which are dual the polytopes
$\Delta_A$ and $\Delta_B$. Moreover, if we consider the lines through
the line segments which are dual to $\conv(m_a,m_b,m_c,m_d)$ and
$\conv(m_a,m_b,m_c,m_f)$ respectively, then these will intersect in a
point $C$ which is dual to the polytope $\conv(m_a,m_b,m_c,m_d,m_f)$
which is not part of the subdivision. Anyway, if we assign to the
points $A$, $B$, and $C$ as weights the lattice area of the
corresponding triangle in $\pi(\Delta)$, e.g.\ $B$ gets as weight the
lattice area $\alpha''$ of $\conv\big((0,0),(0,1),(\alpha'',1)\big)$,
and if we moreover consider the weight of $C$ negatively, since $C$
lies outside the plane segment, then the tropicalisation of the
singular point is the weighted sum of $A$, $B$, and $C$. In the normal
form a simple computation gives $A=(-u,u,0)$, $B=(0,u,0)$ and
$C=(-u,(1+\alpha'')\cdot u,0)$, and
\begin{math}
  \frac{A+\alpha''\cdot B-C}{3}=(0,0,0).
\end{math}
We could thus interpret the tropicalisation of the singular point as a
\emph{virtual weighted barycenter}  of the virtual triangle
$ABC$. \label{page:virtualbarycenter}
\begin{figure}[h]
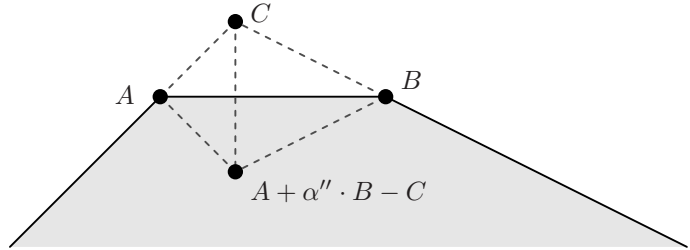

  \centering
  \begin{texdraw}
    \drawdim cm
    \setgray 1
    \move (-3 -3) \lvec (-1 -1) \lvec (2 -1) \lvec (6 -3) \lfill f:0.9
    \setgray 0
    \move (-3 -3) \lvec (-1 -1) \lvec (2 -1) \lvec (6 -3)
    \setgray 0.3
    \lpatt (0.067 0.1)
    \move (-1 -1) \lvec (0 0) \lvec (2 -1)
    \move (-1 -1) \lvec (0 -2) \lvec(2 -1)
    \move (0 0) \lvec (0 -2)
    \move (-1 -1) \fcir f:0 r:0.1
    \move (0 0) \fcir f:0 r:0.1
    \move (2 -1) \fcir f:0 r:0.1
    \move (0 -2) \fcir f:0 r:0.1
    \htext (-1.6 -1.1){$A$}
    \htext (2.2 -0.9){$B$}
    \htext (0.2 0){$C$}
    \htext (0.2 -2.4){$A+\alpha''\cdot B-C$}
  \end{texdraw}
  \caption{The singular point at the virtual barycenter.}
  \label{fig:baricenterE2}
\end{figure}

In our classification we need not consider the Case (b) in
\ref{prop:classification(E)2}, since there the weight class $C$ in
$\Trop(\Ker(A))$ corresponding to this situation is defective
because $\Span(C)$ intersects the lineality space in the vector
corresponding
to the $y$-coordinates of the point configuration.

The Case (c) in Proposition \ref{prop:classification(E)2}
differs from Case (a) by the fact that the points $A$, $B$, and
$C$ all coincide, and that the plane segment corresponding to the
circuit has only one vertex. However, it remains true that the
tropicalisation of the singular point is the weighted sum of $A$, $B$,
and $C$.


Finally, if the Newton polytope contains further points
the situation becomes more complicated. The polytopes $\Delta_A$ and
$\Delta_B$ might be subdivided further, and consequently the vertices
$A$ and $B$ might be cut off, similar to the situation described in
Figure \ref{fig:baricenter2}. As in Subsection \ref{subsec:E1} we can
still identify the virtual points $A$, $B$, and $C$ and their weighted
sum is the tropicalisation of the singular point.

\begin{example}\label{ex:barycenter2}
  A concrete example for this behaviour is the singular point
  $G=(0,0,0)$ on the tropical surface in Example
  \ref{ex-thomas}. Here
  \begin{displaymath}
    \begin{array}{lll}
      m_a=(0,0,0),& m_b=(0,0,1),& m_c=(0,0,2),  \\
      m_d=(1,2,1),&  m_e=(0,1,0), & m_f=(1,0,0),
    \end{array}
  \end{displaymath}
  and one further point $m_g=(-1,-1,0)$. Note that the points
  $m_d,m_e,m_f$ are all on the same side of the plane $x+y=0$ through
  the circuit.  The circuit $m_a,m_b,m_c$
  corresponds then to quadrangle $ABCD$ (see Figure \ref{fig:w-barycenter1},
  \begin{figure}[h]
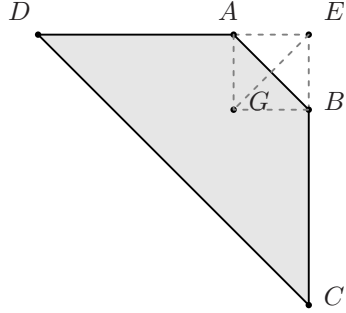

    \centering
    \begin{texdraw}
      \drawdim cm \relunitscale 0.2
      \move (5 -13) \lvec (-13 5) \lvec (0 5) \lvec (5 0)
      \lvec (5 -13) \lfill f:0.9
      \lpatt (0.3 0.5)
      \setgray 0.5
      \move (5 0) \lvec (5 5) \lvec (0 5)
      \move (5 0) \fcir f:0 r:0.2
      \move (0 5) \fcir f:0 r:0.2
      \move (5 5) \fcir f:0 r:0.2
      \move (-13 5) \fcir f:0 r:0.2
      \move (5 -13) \fcir f:0 r:0.2
      \move (0 0) \fcir f:0 r:0.2
      \htext (-15 6) {$D$}
      \htext (6 -13) {$C$}
      \htext (-1 6) {$A$}
      \htext (6 0) {$B$}
      \htext (6 6) {$E$}
      \htext (1 0) {$G$}
      \move (0 5) \lvec (0 0) \lvec (5 -0)
      \move (5 5) \lvec (0 0)
    \end{texdraw}
    \caption{The singular point as barycenter.}
    \label{fig:w-barycenter1}
  \end{figure}
  where the vertices $A=(0,5,0)$ and $B=(5,0,0)$ correspond to the
  polytopes $\Delta_A=\conv\{m_a,m_c,m_d,m_e\}$ respectively
  $\Delta_B=\conv\{m_a,m_c,m_d,m_f\}$ in the subdivision. The polytope
  $\Delta_E=\conv\{m_a,m_c,m_e,m_f\}$, however, is not part of the
  subdivision and defines only a virtual point $E=(5,5,0)$, which is
  the intersection of the two lines determined by the facets
  $\conv\{m_a,m_c,m_e\}$ and $\conv\{m_a,m_c,m_f\}$ of $\Delta_A$
  respectively $\Delta_B$. In this situation and
  \begin{displaymath}
    G=\frac{1}{3}\cdot (C+D-E)
  \end{displaymath}
  is the virtual weighted barycenter of this virtual triangle in the tropical
  surface. Note here, that the virtual vertex $E$ comes with a
  negative weight since it lies outside the plane segment dual to the
  circuit even if we only consider the points $m_a,\ldots,m_f$. Note
  also, that the plane segment dual to the circuit is bounded due to
  the presence of the additional point $m_g$.
\end{example}

\subsection{Weight class as in Lemma \ref{lem-chains}(d), circuit (E) of Remark \ref{rem-circuit}}\label{subsec:circuitE}

\begin{figure}[h]
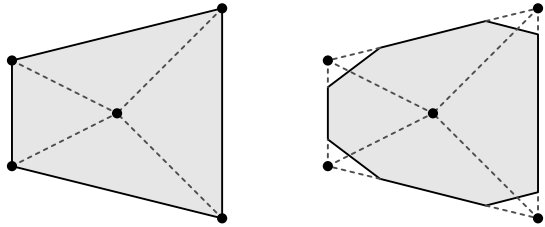

  \centering
  \begin{texdraw}
    \drawdim cm \relunitscale 0.7
    \setgray 1
    \move (-2 -1) \lvec (-2 1) \lvec (2 2) \lvec (2 -2) \lvec (-2 -1) \lfill f:0.9
    \setgray 0
    \move (-2 -1) \lvec (-2 1) \lvec (2 2) \lvec (2 -2) \lvec (-2 -1)
    \setgray 0.3
    \lpatt (0.067 0.1)
    \move (-2 -1) \lvec (0 0) \lvec (-2 1)
    \move (2 2) \lvec (0 0) \lvec(2 -2)
    \move (-2 -1) \fcir f:0 r:0.1
    \move (0 0) \fcir f:0 r:0.1
    \move (-2 1) \fcir f:0 r:0.1
    \move (2 2) \fcir f:0 r:0.1
    \move (2 -2) \fcir f:0 r:0.1
  \end{texdraw}
  \hspace{1cm}
  \begin{texdraw}
    \drawdim cm \relunitscale 0.7
    \setgray 1
    \move (-2 -0.5) \lvec (-2 0.5) \lvec (-1 1.25) \lvec (1 1.75)
    \lvec (2 1.5) \lvec (2 -1.5) \lvec (1 -1.75) \lvec (-1 -1.25) \lvec (-2 -0.5) \lfill f:0.9
    \setgray 0
    \move (-2 -0.5) \lvec (-2 0.5) \lvec (-1 1.25) \lvec (1 1.75)
    \lvec (2 1.5) \lvec (2 -1.5) \lvec (1 -1.75) \lvec (-1 -1.25) \lvec (-2 -0.5) \lfill f:0.9
    \setgray 0.3
    \lpatt (0.067 0.1)
    \move (-1 -1.25) \lvec (-2 -1) \lvec (-2 -0.5)
    \move (-1 1.25) \lvec (-2 1) \lvec (-2 0.5)
    \move (1 1.75) \lvec (2 2) \lvec (2 1.5)
    \move (1 -1.75) \lvec (2 -2) \lvec (2 -1.5)
    \move (-2 -1) \lvec (0 0) \lvec (-2 1)
    \move (2 2) \lvec (0 0) \lvec(2 -2)
    \move (-2 -1) \fcir f:0 r:0.1
    \move (0 0) \fcir f:0 r:0.1
    \move (-2 1) \fcir f:0 r:0.1
    \move (2 2) \fcir f:0 r:0.1
    \move (2 -2) \fcir f:0 r:0.1
  \end{texdraw}
  \caption{The trapeze with the singular point as its midpoint, and the more general situation.}
  \label{fig:trapez}
\end{figure}

Let $F_{s-4}'=\{a,b,c\}$, $F_j'=\{d,e\}$ and $F_i'=\{f,g\}$. We
assume without restriction that $m_a=(0,0,0)$, $m_b=(0,0,1)$
and
$m_c=(0,0,2)$. Dual to this circuit is then as before a
$2$-dimensional polyhedron satisfying $z=0$. We know that in this
situation, the points $m_d$ and $m_e$ lie in a plane with the line
$\{x=y=0\}$, we can assume that this plane satisfies $y=0$. Let us
first assume that $m_d$ and $m_e$ lie on different sides of the
line, i.e. we assume that $m_d$ has positive $x$-coordinate and
$m_e$ has negative $x$-coordinate. Then the triangle with vertices
$m_a$, $m_c$ and $m_d$ (resp.\ $m_e$) will be a face of a polytope
in the subdivision. If $m_d$ or $m_e$ had integral distance bigger
one from the circuit, this face would contain extra lattice points,
contradicting our assumption that the surface is of
maximal-dimensional geometric type. It follows that $m_d$ has $x$-coordinate
$1$ and $m_e$ has $x$-coordinate $-1$. Also, the triangle spanned by
$m_a$, $m_c$ and $m_f$ (resp.\ $m_g$) are faces of the subdivision
and thus $m_f$ and $m_g$ must have integral distance one to the
plane $\{y=0\}$. Let us first assume $m_f$ has $y$-coordinate $1$
and $m_g$ has $y$-coordinate $-1$. Assume first that the subdivision
locally contains only the polytopes $\conv(m_a,m_c,m_d,m_f)$,
$\conv(m_a,m_c,m_e,m_f)$, $\conv(m_a,m_c,m_d,m_g)$ and
$\conv(m_a,m_c,m_e,m_g)$. Then corresponding to this part of the
subdivision we have a quadrangle on the surface. Let us solve for
the $(x,y)$-coordinates of the four vertices. Assume
$m_d=(1,0,\gamma)$, $m_e=(-1,0,\gamma')$, $m_f=(\alpha,1,\gamma'')$
and $m_g=(\alpha',-1,\gamma''')$. Let us denote by
$u=u_{m_a}-u_{m_d}$ the difference of the weights of $m_a$ and $m_d$
and by $w=u_{m_d}-u_{m_f}$ the difference of the weights of $m_d$
and $m_f$. Then the coordinates of the four vertices are
$A=(u,w+(1-\alpha)u)$, $B=(-u,w+(1+\alpha)u)$,
$C=(u,-w+(\alpha'-1)u)$ and $D=(-u,-w-(1+\alpha')u)$. That is, the
quadrangle is a trapeze with the singular point
$(0,0,0)=\frac{A+B+C+D}{4}$ as its midpoint, as depicted in Figure
\ref{fig:trapez} on the left.

If the subdivision contains more polytopes than just these four
locally around the circuit, then we get a polygon with
more sides.
The four edges of the trapeze are still present, and the singular
point is still the midpoint. We can thus extent the cell to a virtual trapeze cell. This more general situation is depicted
in Figure \ref{fig:trapez} on the right.

If $m_d$ and $m_e$ are on the same side of the circuit in the plane
$\{y=0\}$,
then they must both be of integral distance one, and they
form a quadrangle with the circuit which is a face of the
subdivision. Thus the dual subdivision does not correspond to a cone
of the secondary fan of codimension $1$, and we do not consider the
situation.
Analogously, if $m_f$ and $m_g$ are on the same side of the plane
$\{y=0\}$,
they must both have integral distance one to $\{y=0\}$.
However, since the edge connecting $m_f$ and $m_g$ and the circuit
do not need to lie in a plane, it may be that only one of the points
$m_f$ or $m_g$ forms a facet of the subdivision with the circuit. In
this case, the dual subdivision corresponds to a cone of codimension
$1$. However, since the span of the corresponding weight class
intersects the rowspace of $A$ non-trivially (both contain the
vector of $x$-coordinates of the points $m\in \mathcal{A}$), this
weight class is defective and we do not consider the situation.


\end{document}